\documentclass[11pt]{amsart}
\usepackage{amsmath, amscd}

\numberwithin{equation}{section}

\theoremstyle{plain}

\newtheorem{Prop}{Proposition}[section]
\newtheorem{Thm}[Prop]{Theorem}		
\newtheorem{Cor}[Prop]{Corollary}
\newtheorem{Lem}[Prop]{Lemma}
\newtheorem{Def}[Prop]{Definition}
\newtheorem{Rem}[Prop]{Remark}
\newtheorem{Claim}{Claim}

\newcommand{\CBbb}{\mathbb C}
\newcommand{\RBbb}{\mathbb R}
\newcommand{\ZBbb}{\mathbb Z}

\newcommand{\HBbb}{\mathbb H}

\newcommand{\Iso}{\mathop{\rm Iso}\nolimits}

\newcommand{\M}{{\bf M}_{g,n}}
\newcommand{\Mbar}{\overline{{\bf M}}_{g,n}}
\newcommand{\Mprime}{{\bf M}_{g,n}^\prime}
\newcommand{\Mbarprime}{\overline{{\bf M}}_{g,n}^\prime}

\newcommand{\T}{{\bf T}_{g,n}}
\newcommand{\Tbar}{{\overline{\bf T}_{g,n}}}
\newcommand{\Tbarprime}{{\widehat{\bf T}}}
\newcommand{\mcg}{{\Gamma}_{g,n}}

\newcommand{\scc}{\mathcal{C}(\Sigma_{g,n})}
\newcommand{\scck}{\mathcal{C}_k(\Sigma_{g,n})}
\newcommand{\scckprime}{\mathcal{C}_{k'}(\Sigma_{g,n})}
\newcommand{\ck}{{\bf c}_{(k)}}
\newcommand{\ckprime}{{\bf c}^{\prime}_{(k')}}
\newcommand{\MF}{\mathcal{M}\mathcal{F}}
\newcommand{\PMF}{\mathcal{P}\mathcal{M}\mathcal{F}}
\newcommand{\Fmin}{\mathcal{F}_{\min}}
\newcommand{\PFmin}{\mathcal{P}\mathcal{F}_{\min}}
\newcommand{\WF}{{Top}[F]}
\newcommand{\ZF}{{Gr}[F]}
\newcommand{\dt}{\mathop{{\rm d}_{\rm T}}\nolimits}
\newcommand{\dwp}{\mathop{{\rm d}_{\rm WP}}\nolimits}
\newcommand{\twp}{\mathop{{\mathfrak L}_{\rm WP}}\nolimits}
\newcommand{\distortion}{\mathop{\delta_{H,{\mathcal G}}}\nolimits}
\newcommand{\distwp}{\mathop{{\rm dist}_{\rm WP}}\nolimits}
\newcommand{\graft}{ds^2_{t,{\rm graft}}}

\begin{document}



\title[Classification of Weil-Petersson Isometries]
	{Classification of Weil-Petersson Isometries}

\author[Daskalopoulos]{ Georgios Daskalopoulos}

\address{Department of Mathematics \\
		Brown University \\
		Providence,  RI  02912}

\thanks{G.D. supported in part by NSF grant DMS-9803606}

\email{daskal@math.brown.edu}

\author[Wentworth]{ Richard Wentworth}

\address{Department of Mathematics \\
   Johns Hopkins University \\
   Baltimore, MD 21218}

\thanks{R.W. supported in part by NSF grant DMS-9971860}

\email{wentworth@jhu.edu}


\begin{abstract} 
This paper contains two main results.  The first is the existence of an equivariant Weil-Petersson geodesic in
Teichm\"uller space for any choice of pseudo-Anosov mapping class.  As a consequence one obtains a classification of
the elements of the mapping class group as Weil-Petersson isometries which is parallel to the Thurston classification. 
The second result concerns the asymptotic behavior of these geodesics.  It is shown that  geodesics that are
equivariant with respect to independent pseudo-Anosov's diverge.  It follows that subgroups of the mapping class group
which contain independent pseudo-Anosov's act in a reductive manner with respect to the Weil-Petersson geometry.  This
implies an existence theorem for equivariant harmonic maps to the metric completion.
\end{abstract}

$\hbox{}$


\maketitle

\thispagestyle{empty}


\baselineskip=16pt

\section{Introduction}

The Thurston classification of surface diffeomorphisms has strong analogies with
the classification of isometries of  complete nonpositively curved manifolds.   If one regards mapping classes as
isometries with respect to the Teichm\"uller metric then this point of view is made precise by Bers' theorem on the
existence of  Teichm\"uller geodesics that are equivariant with respect to irreducible, nonperiodic diffeomorphisms \cite{B}.  However, the
Teichm\"uller metric is known not to have nonpositive curvature \cite{M1}, so the analogy is not quite germane.

The Weil-Petersson metric does have negative curvature, but it is not complete.  Nevertheless, Wolpert showed that  any two points can be joined
by a Weil-Petersson geodesic \cite{W4}, and it then follows on general principles that  the metric completion of Teichm\"uller space is a
nonpositively curved length space in the sense of Alexandrov (an \emph{NPC} or $CAT(0)$ space).  It is therefore natural and desirable
to recover the classification of mapping classes from the perspective of Weil-Petersson isometries, and that is the first goal of this paper. Note that  by
a recent result of Masur and Wolf, the mapping class group   essentially accounts for all Weil-Petersson isometries
\cite{MW}.  

Let $\mcg$ denote
the mapping class group of a  closed, compact, oriented surface $\Sigma_g$ of genus $g$  with $n$ marked
points. Let $\T$ denote the Teichm\"uller space of $\Sigma_g$ with its marked points, and equip $\T$ with the Weil-Petersson metric.
Then $\mcg$ acts isometrically and properly discontinuously on $\T$ with quotient $\M$, the
Riemann moduli space of curves.  We will denote the metric completion of $\T$ by $\Tbar$.  The extended action of $\mcg$ on $\Tbar$ has quotient
$\Mbar$, the Deligne-Mumford compactification of $\M$.   The first main result of this paper is the following:

\begin{Thm}  \label{T:equivariantgeodesic}
If $\gamma\in\mcg$ is pseudo-Anosov there is a unique $\gamma$-equivariant complete Weil-Petersson geodesic in $\T$.
\end{Thm}

\noindent  By analogy with hyperbolic geometry, we will refer to this geodesic as the {\bf axis} of $\gamma$, and we will denote it by $A_\gamma$. 
The projection of $A_\gamma$ to $\M$ is then a closed geodesic representative of the class of $\gamma$ in the (orbifold) fundamental group of $\M$.

Let us briefly outline the
proof of Theorem \ref{T:equivariantgeodesic}: the first step is a preliminary statement about geodesics in $\Tbar$.  More precisely, we prove in
Proposition \ref{P:wpgeodesic} that a $\gamma$-equivariant geodesic in $\Tbar$ must actually lie in $\T$, provided $\gamma$ is an irreducible mapping
class.  This result is a combination of two ideas.  The first, Theorem \ref{T:interior}, is that any geodesic from a point in $\T$ to the boundary $\partial
\T$ intersects the boundary only at its endpoint.  The second, Theorem \ref{T:corners}, is that the different strata of the boundary intersect transversely.
For example, a geodesic between different components of the top dimensional boundary stratum of $\Tbar$ must pass through the interior $\T$
(see also Remark \ref{R:example}).  The details are given in Section
\ref{S:wp}.

We then consider a length minimizing sequence of equivariant paths $u_j:\RBbb\to\T$ with uniform modulus of
continuity, and let $v_j$ denote the quotients maps to $\M$.  After passing to a subsequence, the $v_j$ converge to a limiting path
$v_\infty:\RBbb\to\Mbar$.  This follows from the compactness of $\Mbar$.  It is an important fact that the space  $\Tbar$ is not even \emph{locally} compact
near points on the boundary
$\partial\T$, so extracting a convergent subsequence directly in $\T$ (or $\Tbar$) is difficult. The second important step in the proof is to show that
$v_\infty$ is the quotient of an equivariant path
$\tilde v_\infty$ to
$\Tbar$.  This is rather technical, and the proof occupies all of Section \ref{S:existence}.  The idea is to inductively lift $v_\infty$ from the lower
dimensional strata $\varepsilon$-close to the path $u_j$.  The parameter $\varepsilon$ is chosen sufficiently small with respect to the injectivity
radius in a compact piece of the lower dimensional stratum so that continuity of the lift from stratum to stratum is preserved. The argument also
shows that the lift is equivariant.   The initial sequence can be taken to be piecewise geodesic, and the geodesic convexity of $\T$ is used
repeatedly.  

Theorem \ref{T:equivariantgeodesic} leads to the desired classification via the {\bf translation
length} $\twp(\gamma)$ of a mapping class $\gamma\in\mcg$, regarded as an isometry of $\T$ with respect to the Weil-Petersson metric (see
(\ref{E:translationlength})).
We refer to
\cite{BGS} for a discussion of how this works for general NPC manifolds and to \cite{Ab} for the case of the Teichm\"uller metric.
 As in those cases,
{\bf semisimple} isometries, i.e.\ those for which the translation length is attained at some point, play an important role. The four possibilities
-- semisimple or not, $\twp$ equal to zero or not -- fit in conveniently with the Thurston classification of surface diffeomorphisms.  According to
Thurston, a mapping class  is either {\bf periodic},  {\bf reducible}, or {\bf pseudo-Anosov} (see Section
\ref{S:mcg} below).  Those  reducible diffeomorphisms that are
 periodic on the component pieces are called {\bf pseudoperiodic}.
 With this terminology, the classification may be
summarized as in  Table 1
 (see Theorem
\ref{T:classification} for the proof).  
\begin{table}[h] 
\centering
\begin{tabular} {|c||c|c|} \hline
&&\\
 &  {\rm semisimple} &  {\rm not semisimple} \\[.1in] \hline\hline
&&\\ 
$\twp=0$  &  periodic &   strictly pseudoperiodic \\[.1in] \hline 
&&\\ 
$\twp\neq 0$  &    pseudo-Anosov &   reducible but not pseudoperiodic \\ 
&&\\ \hline
\end{tabular}
\bigskip
\caption{\emph{Classification of Weil-Petersson Isometries.}}
\end{table}

 The second goal of this paper is to continue the study, initiated in \cite{DKW} of harmonic maps to Teichm\"uller space.
The equivariant harmonic map problem to complete nonpositively curved manifolds (and metric spaces) has been analyzed in great detail (cf.\ \cite{DO,
D1, C, La, JY, KS2}).  Unlike the Dirichlet problem, or the problem of minimizing in a given homotopy class of maps between compact manifolds,  some
condition is generally required of an isometric action on the target manifold in order to guarantee the existence of energy minimizing equivariant
maps.  There are somewhat different notions in various papers; in the case of a symmetric space target, for example, the action should be
\emph{reductive}.

When it comes to Teichm\"uller space,
the salient competing properties of the Weil-Petersson metric  are its negative curvature and its noncompleteness.  Taken together, these two facts
are a reflection of important subtleties in the structure of the mapping class group.  Nevertheless, in light of Theorem \ref{T:equivariantgeodesic}, and
more generally, of the action of
$\mcg$ on the space of measured foliations, the condition that the action contain two independent pseudo-Anosovs (i.e.\  the image subgroup is
{\bf sufficiently large}, see Definition \ref{D:sufficientlylarge})
emerges as a natural candidate for the analogue of the reductive hypothesis just mentioned (cf.\ \cite{DKW}).   

The argument justifying this point of view depends on an understanding of the 
asymptotic geometry of Weil-Petersson geodesics. More precisely, one needs to show that the axes for independent pseudo-Anosov's
diverge, in an appropriate sense, with respect to the Weil-Petersson distance (see Theorem \ref{T:diverge}).   This involves  Thurston's compactification of
$\T$ by projective measured foliations in a crucial way.

  To state the second main result, we recall a key notion in the theory of equivariant harmonic maps (cf.\ \cite[\S 2]{KS2}):
A finitely generated subgroup of
$\mcg$  is called {\bf proper}  if there is a set of generators $\gamma_1,\ldots,\gamma_k$ of the subgroup such that the sublevel sets of the
displacement function 
$$
\delta(\sigma)=\max\left\{ \dwp(\sigma,\gamma_i\sigma) : i=1,\ldots, k\right\}
$$
are bounded (see Definition \ref{D:proper}).  
We apply the same
notion to homomorphisms of groups into $\mcg$ according to the image subgroup. In Section \ref{S:proper} we prove

\begin{Thm} \label{T:proper}
For finitely generated subgroups of $\mcg$, sufficiently large $\Longrightarrow$ proper.
\end{Thm}

This result is important for the harmonic map problem. In the context of the metric completion, the Sobolev theory of
finite  energy   maps from Riemannian domains  to $\Tbar$ has been developed by Korevaar and Schoen \cite{KS1} (see also Jost
\cite{J}).  The important result here  is that energy minimizing (harmonic) maps for both the Dirichlet and equivariant problems are Lipschitz
 continuous.  
The following is an immediate consequence of Theorem \ref{T:proper} and \cite[Theorem 2.1.3 and Remark 2.1.5]{KS2}:

\begin{Cor} \label{C:harmonicmap}
Let $M$ be  a finite volume complete Riemannian manifold with universal
cover $\widetilde M$ and $\rho:\pi_1(M)\to \mcg$ a homomorphism.  Assume that $\pi_1(M)$ is finitely generated and that there exists a
$\rho$-equivariant map $\widetilde M\to\Tbar$  with finite energy (i.e.\ as a map $M\to\Mbar$). Then if $\rho$ is sufficiently large there exists a
finite energy
$\rho$-equivariant harmonic map $u:\widetilde M\to\Tbar$.
\end{Cor}

The sufficiently large condition is flexible enough to provide a wide range of examples so that harmonic maps may be a useful tool in the
study of homomorphisms into mapping class groups.  For instance, by  \cite[Theorem 4.6]{McP}
any irreducible subgroup of $\mcg$ that is not sufficiently large is either finite or virtually cyclic.  It then follows from Corollary
\ref{C:harmonicmap}
that if $M$ is a compact Riemannian manifold and $\rho:\pi_1(M)\to \mcg$ is any homomorphism,
 there
exists a $\rho'$-equivariant harmonic map
$\widetilde M\to \Tbar$, where $\rho'$ is the restriction of $\rho$ to some finite index subgroup of $\pi_1(M)$. 
Possible applications along these lines  are  proofs of certain rigidity results for homomorphisms 
of lattices in Lie groups to mapping class groups and existence results for harmonic representatives of the classifying maps associated to symplectic
Lefschetz pencils. In order to carry out this program, a regularity result for harmonic maps to $\Tbar$ is required.  This issue will be considered in
a future paper.

\smallskip
\noindent
\emph{Acknowledgements.} Theorems \ref{T:equivariantgeodesic} and \ref{T:proper} were first stated by Sumio Yamada in \cite{Y1,Y2}.
We have been strongly influenced by his point of view, but the proof in \cite{Y1} is incomplete, and we have been unable to
carry through his approach.  The methods used in this paper are therefore different and independent from those of Yamada.  We are grateful to
Scott Wolpert for  many helpful suggestions and  for comments on earlier drafts of this paper, and to the referee for a
 careful reading of the manuscript.\footnote{After this paper was submitted for publication, we
received the preprint \cite{W6} which provides an alternative discussion of these results.} 


\section{Mapping Class Groups and Measured Foliations}       \label{S:mcg}


In this section, we review some basic facts about measured foliations and mapping class groups that will be required throughout the paper.  A good
reference for much of this material is \cite{FLP}.   As in the Introduction,  $\Sigma_g$ denotes a closed, compact, oriented surface of genus $g$.  We
will denote a choice of distinct points
$\{p_1,\ldots,p_n\}\subset\Sigma_g$ by $\Sigma_{g,n}$.  The {\bf mapping class group} $\mcg$ is the group of isotopy classes of orientation
preserving homeomorphisms of $\Sigma_{g,n}$ that permute the points $\{p_1,\ldots,p_n\}$.  The equivalence is through isotopies 
fixing the points. We will usually assume that $3g-3+n\geq 0$.

Let $\scc$ denote the set of isotopy classes of simple closed essential nonperipheral curves on $\Sigma_{g,n}$.  In addition, for each integer $k\geq
1$, we let
$\scck$ denote the set of collections ${\bf c}_{(k)}=\{c_1,\ldots,c_k\}$ of $k$ distinct isotopy
classes of  simple closed  curves in $\Sigma_{g,n}$ that are each essential and nonperipheral and which have representatives that are mutually
disjoint.  Given $c$ and $c'$ in $\scc$, let
 $i(c,c')$ denote the geometric
intersection number between them, and for $c\in\scc$ and $\sigma\in\T$,  let $\ell_c(\sigma)$ denote the length of a geodesic representative of $c$
with respect to the complete hyperbolic metric on $\Sigma_{g,n}$ in the class $\sigma$.

Let $\MF$ denote the space of measured foliations on $\Sigma_{g,n}$ modulo isotopy and Whitehead
equivalence.  By an abuse of notation, we shall typically denote a measured foliation and its underlying foliation by the same letter, e.g.\ $F\in\MF$,
assuming that the transverse measure on $F$ is understood.
We do not consider the ``zero measure" to be an
element of  $\MF$.
We can scale the measure on a foliation $F$  by a real number $r>0$ to obtain a new measured foliation, which we denote by $rF$.  We say that $F$ and 
$F'$ are {\bf projectively
equivalent} if $F'=rF$ for some $r$.  The projective equivalence class of a measured foliation $F$ will be denoted by $[F]$, and  the set of
all projective equivalence classes will be denoted by
$\PMF$.

For $F\in\MF$
and
$c\in\scc$, we have an intersection number $i(F,c)$, which by definition is the infimum of the total measures of
representatives of $c$ that are quasi-transverse to $F$.  
Unbounded sequences $\sigma_j\in \T$ have subsequences converging  in $\PMF$ in the
following sense:  there are positive numbers $r_j\to 0$ and $F\in\MF$ such that for each
$c\in\scc$,
$r_j\ell_c(\sigma_j)\to i(F,c)$. Any other sequence $r_j^\prime$ for which
$r_j^\prime\ell_c(\sigma_j)$ converges on all $c$ (and not identically to zero)
gives rise to a measured foliation that is projectively equivalent to $F$.  We shall
refer to this convergence as the {\bf Thurston topology} and denote it by
$\sigma_j\to [F]$. When we want to emphasize convergence in $\MF$, we write $r_j\sigma_j\to F$.

 For any $c\in\scc$ and  $r>0$, there is an associated $rc\in\MF$ such that
$i(rc, c')=ri(c, c')$ for all $c'\in\scc$.
 An important fact is that the intersection number for curves and between curves and
measured foliations extends continuously to a nonnegative function $i(\cdot,\cdot)$ on
$\MF\times \MF$ such that
$i(rF,r'F')=rr'i(F,F')$.  By analogy with curves, a pair $F, F'\in\MF$ is called {\bf transverse} if $i(F,F')\neq 0$.
A measured foliation $F$ is called {\bf minimal} if $i(F,c)>0$ for every
$c\in\scc$.  Let $\Fmin\subset\MF$ denote the subset  of minimal
foliations.  Since the condition is independent of the projective factor, we also obtain a subset $\PFmin\subset\PMF$ of projective minimal
foliations.   For any
$[F]\in\PMF$, we let 
\begin{align}
\WF&=\left\{ [G]\in\PMF : i(F,G)=0\right\} \label{E:top} \\
\ZF &= \left\{ [G]\in\PMF : i(G,c)=0  \iff i(F,c)=0\; \text{for all}\; c\in\scc
\right\} \label{E:gr}
\end{align}
Note that these definitions are independent of the choice of representatives of the classes $[F]$ and $[G]$.
 The sets $\ZF$ give a countable partition of
$\PMF\setminus\PFmin$, while for
 $[F]\in\PFmin$, the  condition defining $\WF$  is an equivalence relation.  Indeed, the set $\WF$
consists precisely of the projective equivalence classes of those measured foliations whose underlying foliations are topologically equivalent to that of
$F$, but whose measures may be different.

Thurston's classification of surface diffeomorphisms may be described in terms of the natural action of $\mcg$ on $\MF$ and $\PMF$:
An element $\gamma\in\mcg$ is called {\bf reducible} if $\gamma$ fixes some collection $\ck\in\scck$. 
 It is called
{\bf pseudo-Anosov} if there is  $r>1$ and transverse measured
foliations
$F_+$, $F_-$ on $\Sigma_{g,n}$ such that $\gamma F_+\simeq r F_+$, and $\gamma F_-\simeq r^{-1} F_-$.  
 $F_+$ and $F_-$ are called the {\bf stable} and {\bf unstable} foliations of $\gamma$, respectively. The classification states that
any
$\gamma\in\mcg$ is either periodic (i.e.\ finite order), infinite order and reducible, or pseudo-Anosov.  Moreover, these are mutually exclusive
possibilities. 

 The stable and unstable foliations of a pseudo-Anosov are both minimal and carry a unique transverse
measure, up to projective equivalence.   We will need the following facts:
\begin{align}
& \gamma \text{ pseudo-Anosov and } [F]\not\in\PFmin \Longrightarrow [\gamma F]\not\in\ZF\ .  \label{E:zf} \\
& \gamma  \text{ pseudo-Anosov  and } [F]\in\PFmin\setminus\left\{[F_+],[F_-]\right\} \Longrightarrow
[\gamma F]\not\in\WF\ .  \label{E:wf}
\end{align}

Finally, a pseudo-Anosov element fixes \emph{precisely} two points in $\PMF$; namely the points $[F_+]$ and $[F_-]$ represented by the stable and unstable
foliations.   We say that pseudo-Anosov's  are {\bf independent} if their fixed point sets in $\PMF$ do not coincide.

With this background we are now prepared to make the following important

\begin{Def}[{\cite[p.\ 142]{McP}}]  \label{D:sufficientlylarge} 

A subgroup of
$\mcg$ is  {\bf sufficiently large} if it contains two independent pseudo-Anosov's.  
\end{Def}


\section{The Weil-Petersson Metric}       \label{S:wp}


Let $\T$ denote the Teichm\"uller space of $\Sigma_{g,n}$, and $\M=\T/\mcg$ the Riemann moduli space of curves. The important metric
 on $\T$ for this paper is the
 Weil-Petersson metric, which we denote by $\dwp$.  We will also have occasion to use the  Teichm\"uller metric, which we denote by $\dt$.   
This is a complete Finsler metric with respect to which
$\mcg$ also acts by isometries.

Let us briefly recall the definitions: 
 the cotangent space $T^\ast_\sigma \T$ is identified
with the space of meromorphic quadratic differentials on $(\Sigma_{g,n},\sigma)$ that have at most simple poles at the marked points $\{p_1,\ldots,p_n\}$.
Here we abuse notation slightly, and regard the point $\sigma\in\T$ as giving a conformal structure on $\Sigma_{g,n}$ with its marked points.  The
complete hyperbolic metric on $(\Sigma_{g,n},\sigma)$ can be expressed in local conformal coordinates as $ds^2=\rho(z)|dz|^2$.  Similarly, a
quadratic differential has a local expression $\Phi=\varphi(z)dz^2$.
 Then for $\Phi\in
T^\ast_\sigma \T$, the
{\bf Weil-Petersson cometric} is  given by
\begin{equation}
\Vert\Phi\Vert^2_{WP}=\int_\Sigma |\varphi(z)|^2\rho(z)^{-1} dxdy\ .
\end{equation}
The  {\bf Teichm\"uller cometric} is given by 
\begin{equation}
\Vert\Phi\Vert_T=\int_\Sigma |\varphi(z)|dxdy\ .
\end{equation}

\noindent 
 For a point
$\sigma\in\T$ and a subset $A\subset\T$, we set
$$
\distwp(\sigma,A)=\inf\left\{ \dwp(\sigma,\sigma') : \sigma'\in A\right\} \ .
$$

We record the following important facts:

\begin{Prop}[{cf.\ \cite{Ah,W1,W3,W4,Chu,R,Tr}}]  \label{P:wolpert}
\par\noindent
\begin{enumerate}
\item  The Weil-Petersson metric is a noncomplete K\"ahler metric with negative sectional curvature.
\item  Any two points in $\T$ may be joined by a unique Weil-Petersson geodesic in $\T$.
\item  Any length function $\ell_c$ is strictly convex.
\item  There is a strictly convex exhaustion function on $\T$.
\end{enumerate}
\end{Prop}

Let $\Tbar$ denote the metric completion of $\T$ with respect to the Weil-Petersson metric.   Further notation: let $\partial\T=\Tbar\setminus\T$ and
$\partial
\M=\Mbar\setminus\M$.   The completion has the following  local description (cf.\ \cite{M2}):  $\partial\T$ is a disjoint union
of smooth connected strata $D(\ck)$.  These are described  by choosing a collection $\ck=\{c_1,\ldots, c_k\}$ of disjoint simple closed essential and
nonperipheral curves  on
$\Sigma_{g,n}$ and forming a nodal surface by collapsing each of the $c_j$ to  points.  
Associated to  a nodal surface is another
  Teichm\"uller
space which is by definition the set of equivalence classes of conformal structures on the normalized (possibly
disconnected) surface, with the preimages of the nodes as additional marked points. It is therefore naturally isomorphic to a product of lower
dimensional Teichm\"uller spaces.
The statement is that the stratum
$D(\ck)$ with its induced metric is isometric to this product of lower dimensional Teichm\"uller spaces with the product
Weil-Petersson metric.   

To describe this in more detail we define an incomplete metric space
\begin{equation} \label{E:modelmetric}
\HBbb=\left\{ (\theta,\xi)\in \RBbb^2 : \xi>0\right\} \ ,\qquad
ds^2_{\HBbb}=4d\xi^2+\xi^6d\theta^2\ .
\end{equation}
Let $\overline\HBbb$ denote the metric completion of $\HBbb$, obtained by adding a single point $\partial\HBbb$  corresponding to
the entire real axis $\xi=0$, and denote the distance function on $\overline\HBbb$ by $d_{\overline\HBbb}$. The completion is then an NPC
space which is, however, not locally compact.  We note the following simple:

\begin{Lem}  \label{L:modelgeodesic}
The geodesic $w:[0,1]\to\overline\HBbb$ from $\partial\HBbb$ to any point $(\theta_1,\xi_1)\in\HBbb$  is given by
$w(x)=(\theta_1,x\xi_1)$.
\end{Lem}

The importance of $\overline\HBbb$ is that it is a model for the normal space to the boundary strata.  More precisely,
let $\sigma\in D(\ck)$. Then the curves in $\ck$ divide $\Sigma_{g,n}$ into a disjoint union of surfaces $\Sigma_{g_1,n_1},
\ldots, \Sigma_{g_N,n_N}$, where the boundary components corresponding to the curves are collapsed to additional marked points.
Set 
\begin{equation} \label{E:teich}
\widehat{\bf T}={\bf T}_{g_1,n_1}\times\cdots\times{\bf T}_{g_N,n_N}\ .
\end{equation} 
Then $\sigma$ corresponds to a point  in $\widehat {\bf T}$, which we will also denote by $\sigma$.
There is an open neighborhood $U$ of $\sigma$ in $\Tbar$ that is homeomorphic to an open neighborhood in the 
product
$\overline\HBbb^k\times\widehat{\bf T}
$, where $\sigma$ is mapped to the point $(\partial\HBbb,\ldots,\partial\HBbb,
\sigma)$.  We will call such a $U$ a {\bf model neighborhood} of $\sigma$.

Next, we describe the asymptotic  behavior of the Weil-Petersson metric near the boundary using a model neighborhood $U$.  We choose complex
coordinates
$\tau_{k+1},\ldots, \tau_{3g-3+n}$ for $\widehat{\bf T}$ near $\sigma$.  Let $(\theta_i,\xi_i)$, $i=1,\ldots, k$  denote
the coordinates for each factor of
$\overline \HBbb$ in the model neighborhood.  With respect to these coordinates, $U\cap D(\ck)$ is given by the equations
$\xi_{1}=\cdots=\xi_{k}=0$.  As  $\xi=(\xi_{1},\ldots,\xi_{k})\to 0$, the Weil-Petersson metric has an expansion:

\begin{align}
ds_{\rm WP}^2 &=\sum  \left(G_{i\bar j}+ \sum_\ell O\left(\xi_\ell^4\right)\right)\, d\tau_i\otimes d\bar \tau_j\notag \\
&\quad +\sum O\left(\xi_j^3\right)\times\left[ d\tau_i\otimes d\xi_j\ \text{or}\ d\bar \tau_i\otimes d\xi_j\right]\notag \\
&\quad +\sum O\left(\xi_j^6\right)\times\left[ d\tau_i\otimes d\theta_j\ \text{or}\ d\bar \tau_i\otimes d\theta_j\right]\notag \\
&\quad +\sum\left(B_i+\sum_\ell O\left(\xi_\ell^4\right)\right)\, 4d\xi_i^2\notag \\
&\quad + \sum_{i\neq j}O\left(\xi_i^3\xi_j^3\right)\, d\xi_i\otimes d\xi_j\label{E:wp} \\
&\quad + \sum \left(B_i+\sum_\ell O\left(\xi_\ell^4\right)\right)\xi_i^6\, d\theta_i^2\notag \\
&\quad + \sum_{i\neq j} O(\xi_i^6\xi_j^6)\, d\theta_i\otimes d\theta_j\ ,\notag
\end{align}

\noindent  where $G_{i\bar j}$ denotes the Weil-Petersson metric on $\widehat{\bf T}$, and the coefficients $B_i\neq 0$ are
continuous functions of $\tau_{k+1},\ldots, \tau_{3g-3+n}$   defined in a neighborhood (which we assume is the model neighborhood) of $
\sigma\subset \widehat{\bf T}$.  In particular, the induced metrics on the factors $\overline \HBbb$ are asymptotic to metrics that are uniformly
quasi-isometric to the model  (\ref{E:modelmetric}). This expansion is due to Yamada \cite{Y1}.  For the sake of completeness, we provide a brief
justification in the Appendix below.

\begin{Rem}  \label{R:npc}
As noted in \cite[Remarks in \S 1.3]{MW}, $\Tbar$ is a complete  nonpositively curved length space in the
sense of Alexandrov, i.e.\
$\Tbar$ is an NPC (or equivalently, a CAT(0)) space with an isometric action of $\mcg$ extending the action on $\T$ (cf.\ \cite[Corollary
3.11]{BH}).  Moreover, the strata
$D(\ck)\subset\partial\T$ are totally geodesically embedded.
\end{Rem}

\begin{Rem} \label{R:ks}
Since $\Tbar$ is a length space, between any two points $\sigma_0$, $\sigma_1$ we can find a rectifiable path $u:[0,1]\to\Tbar$
with 
\begin{equation*}
(\ast)\
\begin{cases}
u(0)&=\sigma_0\\
u(1)&=\sigma_1
\end{cases}
\end{equation*}
whose length $L(u)$ coincides with the distance $\dwp(\sigma_0,\sigma_1)$.  We will often regard $u$ as an energy minimizing path
(in the sense of \cite{KS1}) with respect to the boundary conditions $(\ast)$.  We will denote the energy of such a path by
$E(u)$.
\end{Rem}

Before proceeding, let us introduce some convenient notation regarding $\Tbar$:  For a fixed $k\geq 1$, let $\dot\Delta_k$ denote the union of strata
$D(\ck)$ for all possible collections $\ck\in\scck$, and let $\Delta_k$ denote its closure in the completion $\Tbar$.  Then we have
\begin{equation}  \label{E:stratification} 
\partial\T=\bigcup_{k= 1,\ldots, k_{max}}\dot\Delta_k \ ,\quad
\dot\Delta_k = \Delta_k\setminus\Delta_{k+1}\ . 
\end{equation}

\noindent Here, $k_{max}=3g-3+n$ is the maximal number of disjoint simple closed essential and nonperipheral curves on $\Sigma_{g,n}$.   Note that with this
notation,
$\partial\T=\Delta_1$. It will also be useful to set $\dot\Delta_0=\T$.

To state the main result of this section, we clarify some terminology that we have already used in the statement of Theorem
\ref{T:equivariantgeodesic}: Let $\gamma\in\mcg$ and $u:\RBbb\to\T$.  We will say that $u$ is {\bf equivariant} with respect to
$\gamma$ if
$u(x+L)=\gamma u(x)$ for all $x\in\RBbb$ and some fixed $L>0$.  By convention, 
the length of $u$ will then be taken to mean the length of $u$
restricted to
$[0,L]$. Unless otherwise stated, we will parametrize so that $L=1$.
Also, for
$\gamma\in\mcg$, we define the  {\bf Weil-Petersson translation length} of $\gamma$ by 
\begin{equation}  \label{E:translationlength}
\twp(\gamma)=\inf\left\{\dwp(\sigma, \gamma\sigma): \sigma\in\T\right\}\ .
\end{equation}

The following  characterization of  Weil-Petersson geodesics in terms of the completion 
will be used  in the next section:

\begin{Prop}  \label{P:wpgeodesic}
Suppose $\gamma\in\mcg$ is irreducible.  Let $u:\RBbb\to\Tbar$ be a rectifiable $\gamma$-equivariant path, and suppose that
$L(u)\leq\twp(\gamma)$.  Then
$u$ is a Weil-Petersson geodesic and has image in
$\T$.
\end{Prop}

\noindent   The remainder of this section is devoted to the proof of this result.  It will follow from two properties of the geometry of the Weil-Petersson
completion, Theorems \ref{T:interior} and \ref{T:corners}, that were first stated by Yamada \cite{Y1}.  The Harnack type estimate  in \cite{Y1} does not,
however,  appear to hold in the generality needed.  The proofs given below will proceed along a different line of reasoning. Specifically, we use a
rescaling argument. The main result is the statement that geodesics
in $\Tbar$ between points in different boundary components must pass through the interior $\T$.  We formulate this in Theorem
\ref{T:corners} below.
First we show that geodesics from points in Teichm\"uller space to the boundary touch the boundary only at their
endpoints.

\begin{Thm} \label{T:interior}
Let $u:[0,1]\to\Tbar$ be a geodesic with $u(0)=\sigma_0\in\partial\T$ and $u(1)=\sigma_1\in\T$.  Then $u((0,1])\subset\T$. 
Similarly, if $\sigma_0\in\partial D(\ck)$ and $\sigma_1\in D(\ck)$, then $u((0,1])\subset D(\ck)$.
\end{Thm}

\begin{proof}
First consider the case where $\sigma_0\in\partial\T$ and $\sigma_1\in\T$.
Let 
$
\hat x=\sup\left\{x : u(x)\in \partial \T \right\}
$.
We wish to prove that $\hat x=0$.  Suppose not, and set $\sigma_{\hat x}=u(\hat x)$.  Then there is a collection $\ck\in\scck$ such that
$\sigma_{\hat x}\in D(\ck)$.  By changing the endpoints we may assume that the image of $u$ is contained in a model neighborhood of
$\sigma_{\hat x}$.  

\begin{Claim} \label{C:dck}
$u([0,\hat x])\subset D(\ck)$.
\end{Claim}

\begin{proof}  Let $c\in\ck$, and suppose that $\ell_c u(x)\neq 0$ for some $0\leq x<\hat x$.  We may find $\gamma\in\scc$ such that
$i(\gamma, c)\neq 0$, but $i(\gamma, c')=0$ for any $c'\in\ck$, $c'\neq c$.  Let $\tilde u_j:[x,1]\to\T$ be a sequence of geodesics converging
 in $\Tbar$ to $u\bigr|_{[x,1]}$.  By Proposition \ref{P:wolpert} (3),  $f_j=\ell_\gamma \tilde u_j$ is a sequence of convex functions on
$[x,1]$, and since $u(1)\in\T$ and $\ell_c u(x)\neq 0$, $f_j(1)$ and $f_j(x)$ are  uniformly bounded.  Hence, $f_j(\hat x)$ is also uniformly bounded. 
But this contradicts the fact that
$\ell_c\tilde u_j(\hat x)\to 0$ and $i(\gamma,c)\neq 0$.
\end{proof}

  By definition, $u((\hat x,1])\subset\T$.  Thus, after reparametrizing, we have the following situation:
$
u:[-1,1]\to \overline\HBbb^k\times\Tbarprime 
$
is energy minimizing with respect to its boundary conditions.    Furthermore,
\begin{equation} \label{E:corner}
u([-1,0])\subset \left(\partial\HBbb\right)^k\times\Tbarprime \ , \ \text{and}\ \
u((0,1])\subset \HBbb^k\times\Tbarprime \ .
\end{equation}
We shall derive a contradiction to the existence of a geodesic satisfying (\ref{E:corner}).  The strategy is to use a rescaling argument to
compare Weil-Petersson geodesics to geodesics in $\overline\HBbb^k\times\Tbarprime $ -- where we take 
the product of the model metric on $\overline\HBbb^k$ and the Weil-Petersson metric on $\Tbarprime$ --
and thus derive a contradiction to Lemma \ref{L:modelgeodesic}. 
Without loss of generality, we assume that $B_i(\sigma_{\hat x})=1$, $i=1,\ldots, k$, in the asymptotic expression (\ref{E:wp}).

Let $u^{\prime}$ denote the projection of $u$ onto the first factor $\overline\HBbb^k$.  We also choose a sequence of
rescalings $\varepsilon_j\downarrow 0$, and set $u_j(x)=u(\varepsilon_j x)$,
 $u^{\prime}_j(x)=u^{\prime}(\varepsilon_j x)$. Note that each
$u^{\prime}_j$ satisfies:
\begin{equation} \label{E:boundary}
u^{\prime}_j(x)=\left(\partial\HBbb\right)^k \qquad \text{for all}\ x\in[-1,0]\ .
\end{equation}

\noindent    Let $w_j$ denote the
 map $[-1,1]\to\overline\HBbb^k$ with boundary conditions 
\begin{equation} \label{E:boundaryconditions}
w_j(-1)=u^{\prime}_j(-1)=\left(\partial\HBbb\right)^k\ ,\qquad
w_j(1)=u^{\prime}_j(1)\in\HBbb^k\ ,
\end{equation}
which is harmonic with respect to the model metric (see (\ref{E:modelmetric}) and Lemma \ref{L:modelgeodesic}).  Note that since
$w_j$ is a nonconstant energy minimizer, 
$
0\neq E_{\overline\HBbb^k}(w_j)\leq E_{\overline\HBbb^k}(u^{\prime}_j)
$
for all $j$.  Here, $E_{\overline\HBbb^k}(u)$ denotes the energy of a path in $\overline\HBbb^k$ with respect to the product of the model
metrics
$ds^2_{\HBbb}$.

\noindent 
\begin{Claim}  \label{C:lipschitz}
There is a constant $C$ independent of $j$ such that
$$
\max_{x\in[-1,1]} d_{\overline \HBbb^k}(u^{\prime}_j(x),\left(\partial\HBbb\right)^k)\leq C\,
E_{\overline\HBbb^k}^{1/2}(u^{\prime}_j)\ ,
$$
and
$$
\max_{x\in[-1,1]} d_{\overline \HBbb^k}(w_j(x),\left(\partial\HBbb\right)^k)=d_{\overline \HBbb^k}(w_j(1),\left(\partial\HBbb\right)^k)\leq C\,
E_{\overline\HBbb^k}^{1/2}(u^{\prime}_j)\ .
$$
\end{Claim}

\begin{proof}
From the estimate (\ref{E:wp}) the Weil-Petersson metric on $\overline\HBbb^k\times\Tbarprime$ is quasi-isometric to the product of the model
metric on
$\overline\HBbb^k$ and the Weil-Petersson metric on a lower dimensional Teichm\"uller space.  Let ${\mathcal
B}=\left(\partial\HBbb\right)^k\times\Tbarprime$. Then the quasi-isometry implies  that there is a constant $C$ such that for all $j$ and all
$x\in[-1,1]$, 
\begin{equation} \label{E:qi}
C^{-1}\dwp(u_j(x), {\mathcal B})\leq d_{\overline\HBbb^k}(u_j^{\prime}(x),\left(\partial\HBbb\right)^k)\leq C\dwp(u_j(x), {\mathcal B})\
.
\end{equation}
 Since ${\mathcal B}$ is convex, $u_j$ is a geodesic, and $\Tbar$ is NPC, $x\mapsto\dwp(u_j(x),{\mathcal B})$ is a convex function.  Hence,
\begin{equation} \label{E:convex}
\max_{x\in[-1,1]} \dwp(u_j(x),{\mathcal B})= \dwp(u_j(1),{\mathcal B})\ .
\end{equation}
Combining (\ref{E:qi}) with (\ref{E:convex}), we have\footnote{We follow the standard practice of using $C$ to denote a generic constant whose value
may vary from line to line.}
$$\max_{x\in[-1,1]} d_{\overline \HBbb^k}(u^{\prime}_j(x),\left(\partial\HBbb\right)^k)\leq C d_{\overline \HBbb^k}(u^{\prime}_j(1),
\left(\partial\HBbb\right)^k)\leq C\, L(u^{\prime}_j)\leq C\,  E_{\overline\HBbb^k}^{1/2}(u^{\prime}_j)\ ,
$$
which is the first inequality in the claim.
The second inequality in the claim is immediate from the fact that $w_j$ is a geodesic and $E_{\overline\HBbb^k}(w_j)\leq
E_{\overline\HBbb^k}(u^{\prime}_j)$.
\end{proof}

\begin{Claim} \label{C:energyestimate}
There are sequences $\mu_j$ and $\eta_j$ such that $\mu_j\to 1$ and $\eta_j\, E_{\overline\HBbb^k}^{-1}(u^{\prime}_j)\to 0$ as $j\to
\infty$, and 
$$
0\neq E_{\overline\HBbb^k}(w_j)\leq E_{\overline\HBbb^k}(u^{\prime}_j)\leq\mu_j\, E_{\overline\HBbb^k}(w_j)
+\eta_j\ .
$$
\end{Claim}

\begin{proof}
Let $u_j^{\prime\prime}$ denote the projection of $u_j$ onto the $\Tbarprime$ factor of the product $\overline\HBbb^k\times\Tbarprime$.  The
Weil-Petersson energy may be bounded
\begin{equation} \label{E:errorterm}
E_{\mathcal B}(u_j^{\prime\prime})+\delta_j^{-1}\, E_{\overline\HBbb^k}(u^{\prime}_j)-E(u_j^{\prime},u_j^{\prime\prime})\leq
E(u_j)\leq
E_{\mathcal B}(u_j^{\prime\prime})+\delta_j\, E_{\overline\HBbb^k}(u^{\prime}_j)+E(u^{\prime}_j,u_j^{\prime\prime})\ ,
\end{equation}
where $E(u^{\prime}_j,u_j^{\prime\prime})\geq 0$ is a cross term reflecting the error of the Weil-Petersson metric from the product
metric, and $E_{\mathcal B}(u_j^{\prime\prime})$ denotes the energy of $u_j^{\prime\prime}$ with respect to the metric on ${\mathcal
B}=\left(\partial\HBbb\right)^k\times\Tbarprime$ induced from the Weil-Petersson metric by the inclusion in $\Tbar$.  Such an estimate follows
from  (\ref{E:wp}).  The constants
$\delta_j\geq 1$ arise because the induced metric on
$\overline
\HBbb^k$ is only quasi-isometric to the product of (\ref{E:modelmetric}).  Because the functions $B_i$ in (\ref{E:wp}) corresponding
to the curves
$c_1,\ldots, c_k$ are continuous, and we have assumed 
$B_i(\sigma_{\hat x})=1$, we  conclude that $\delta_j\to 1$ as $j\to\infty$.  Consider the map
$\hat u_j$ with coordinates
$( w_j,u_j^{\prime\prime})$.  This has the same boundary conditions as $u_j$.  
Using (\ref{E:errorterm}) applied to $\hat u_j$, we have
\begin{equation} \label{E:firsthalf}
E(\hat u_j)\leq
E_{\mathcal B}(u_j^{\prime\prime})+\delta_j\, E_{\overline\HBbb^k}(w_j)+E(w_j,u_j^{\prime\prime})\ .
\end{equation}
The left hand side of (\ref{E:errorterm}) is:
\begin{equation} \label{E:secondhalf}
E_{\mathcal B}(u_j^{\prime\prime})+\delta_j^{-1}\, E_{\overline\HBbb^k}(u^{\prime}_j)-E(u_j^{\prime},u_j^{\prime\prime})\leq
E(u_j)\ .
\end{equation}
On the other hand, since  $u_j$ is energy minimizing, $E(u_j)\leq E(\hat
u_j)$.  Using this fact and eqs.\ (\ref{E:firsthalf}) and (\ref{E:secondhalf}), we have: 
\begin{equation} \label{E:claim2}
E_{\overline\HBbb^k}(u^{\prime}_j)\leq
\mu_j E_{\overline\HBbb^k}(w_j)+\delta_j\left(E( w_j,u_j^{\prime\prime})+E(u^{\prime}_j,u_j^{\prime\prime})\right)\ ,
\end{equation}
where $\mu_j=\delta^2_j\to 1$.

By the estimate (\ref{E:wp}) the cross terms in the expansion of the Weil-Petersson metric near $\partial\T$ depend on terms involving
$(u_j)^3_{\xi_i}$, $i=1,\ldots, k$.  Moreover, the uniform Lipschitz bound on $u$ implies that  $|(d/dx)
(u_j^{\prime})_{\xi_i}|$ and 
$(u_j^{\prime})^3_{\xi_i}|(d/dx) (u_j^{\prime})_{\theta_i}|$ are bounded.  Similarly, the norm of $(d/dx) u_j^{\prime\prime}$ with respect to the
induced metric is also uniformly bounded. 
 We then have an estimate:
\begin{align*}
E(u^{\prime}_j,u_j^{\prime\prime})&\leq C\left(\max_{x\in[-1,1]} d_{\overline
\HBbb^k}(u^{\prime}_j(x),\left(\partial\HBbb\right)^k)\right)^3\ ,\\ E(w_j,u_j^{\prime\prime})&\leq C\left(\max_{x\in[-1,1]} d_{\overline
\HBbb^k}(w_j(x),\left(\partial\HBbb\right)^k)\right)^3\ ,
\end{align*}
for some constant $C$ independent of $j$.  Applying Claim \ref{C:lipschitz} we have:
$$
E(u^{\prime}_j,u_j^{\prime\prime})+
E(w_j,u_j^{\prime\prime})\leq C\,  E^{3/2}_{\overline\HBbb^k}(u^{\prime}_j)\ .
$$
 Plugging this into (\ref{E:claim2}) proves the claim.
\end{proof}

We now rescale the metric on $\overline \HBbb^k$ by $E^{-1/2}_{\overline\HBbb^k}(u^{\prime}_j)$ to obtain  metrics
$d_{\overline\HBbb^k}^{(j)}$.  We shall denote the energies with respect to these new metrics by $E^{(j)}_{\overline\HBbb^k}$.  Note
that by definition of this rescaling and Claim 2, we have:
\begin{align}
E^{(j)}_{\overline\HBbb^k}(u^{\prime}_j)&=1\quad\text{for all $j$,} \label{E:uenergy} \\
\lim_{j\to\infty}E^{(j)}_{\overline\HBbb^k}(w_j)&= 1\ .\label{E:wenergy}
\end{align}

\begin{Claim} \label{C:poincare}
$\displaystyle \lim_{j\to\infty} d^{(j)}_{\overline\HBbb^k}(u_j^{\prime}(x), w_j(x)) =0$, uniformly for $x\in[-1,1]$.
\end{Claim}

\begin{proof}
Let $\hat w_j(x)$ be the midpoint along the geodesic from $u_j^{\prime}(x)$ to $w_j(x)$  in
$\overline\HBbb^k$ (with the metric $d^{(j)}_{\overline\HBbb^k}$).  By  \cite[proof of Theorem 2.2]{KS1} it follows that the map $x\mapsto \hat w_j(x)$ is an
admissible competitor to
$w_j$ with the same boundary conditions.  Moreover, we have \cite[eq.\ (2.2iv)]{KS1}:
$$
2E^{(j)}_{\overline\HBbb^k}(\hat w_j)\leq  E^{(j)}_{\overline\HBbb^k}(u^{\prime}_j) + E^{(j)}_{\overline\HBbb^k}(w_j)
-\frac{1}{2}\int_{-1}^1 \left|\frac{d}{ds} d^{(j)}_{\overline\HBbb^k}(u_j^{\prime}(s), w_j(s))\right|^2ds\ .
$$
By (\ref{E:uenergy})  and the fact that $w_j$ is a minimizer:
\begin{align*}
2 E^{(j)}_{\overline\HBbb^k}(w_j)\leq 2E^{(j)}_{\overline\HBbb^k}(\hat w_j)&\leq 1+E^{(j)}_{\overline\HBbb^k}(w_j) 
-\frac{1}{2}\int_{-1}^1 \left|\frac{d}{ds} d^{(j)}_{\overline\HBbb^k}(u_j^{\prime}(s), w_j(s))\right|^2ds\\
\frac{1}{2}\int_{-1}^1 \left|\frac{d}{ds} d^{(j)}_{\overline\HBbb^k}(u_j^{\prime}(s), w_j(s))\right|^2ds
&\leq 1-E^{(j)}_{\overline\HBbb^k}(w_j)\ .
\end{align*}
By (\ref{E:wenergy}), it follows that
$$
\lim_{j\to\infty}\int_{-1}^1 \left|\frac{d}{ds} d^{(j)}_{\overline\HBbb^k}(u_j^{\prime}(s), w_j(s))\right|^2ds = 0 \ .
$$
By (\ref{E:boundaryconditions}) we have for any $x\in[-1,1]$,
$$
d^{(j)}_{\overline\HBbb^k}(u_j^{\prime}(x), w_j(x))\leq 
\int_{-1}^x \left|\frac{d}{ds} d^{(j)}_{\overline\HBbb^k}(u_j^{\prime}(s), w_j(s))\right|ds
\leq \sqrt{2}\left\{\int_{-1}^1 \left|\frac{d}{ds} d^{(j)}_{\overline\HBbb^k}(u_j^{\prime}(s), w_j(s))\right|^2ds\right\}^{1/2}\ .
$$
The claim now follows.
\end{proof}

We continue  with the proof of Theorem \ref{T:interior}.  First, note that since the $w_j$ are constant speed geodesics, we may assume that the
pull-back distance functions $w_j^\ast d_{\overline \HBbb^k}^{(j)}$ converge to a nonzero multiple of the euclidean distance on
$[-1,1]$.  In particular, there is a $c\neq 0$ such that for any $x\in[-1,1]$,
$$
d_{\overline \HBbb^k}^{(j)}(w_j(x),w_j(-1))\to c|x+1|\ .
$$
By Claim \ref{C:poincare} it follows that
$$
d_{\overline \HBbb^k}^{(j)}(u^{\prime}_j(x),w_j(-1))\to c|x+1|\ ,
$$
for all $x$.  But by (\ref{E:boundary}) and (\ref{E:boundaryconditions}),
$$
d_{\overline \HBbb^k}^{(j)}(u^{\prime}_j(x),w_j(-1))=d_{\overline \HBbb^k}^{(j)}(u^{\prime}_j(x),\left(\partial\HBbb\right)^k)=0\ ,
$$
for all $x\in [-1,0]$ and all $j$.
This contradiction rules out the existence of a geodesic satisfying (\ref{E:corner}), and therefore completes the proof of the Theorem for
the first case.  Since $D(\ck)\hookrightarrow\Tbar$ is totally geodesic, the proof for the case where
$\sigma_0\in\partial D(\ck)$ and
$\sigma_1\in D(\ck)$ is identical.
\end{proof}

Next we show that the different strata of the boundary intersect transversely:

\begin{Thm} \label{T:corners}
Let $u:[0,1]\to\Tbar$ be a geodesic.  Suppose $u$ is contained in some $\Delta_k$, but is not completely contained in
$\Delta_{k+1}$.  Then $u$ is contained in the closure of a single connected component of $\dot\Delta_k$.
\end{Thm}

\begin{Rem} \label{R:example}
To clarify the statement of this result, let us give an example.  Consider disjoint nonisotopic simple closed essential curves $\{c_1, c_2\}$ on
a closed compact surface $\Sigma_g$.  Let
$\sigma_1$ denote a point in the boundary component $D(c_1)$ of $\overline{\bf T}_g$ corresponding to pinching $c_1$.  Similarly, let $\sigma_2$ denote a
point in the boundary component $D(c_2)$ of $\overline{\bf T}_g$ corresponding to pinching $c_2$.   Since $c_1$ and $c_2$ are disjoint, the
intersection of the closures
$\overline D(c_1)\cap\overline D(c_2)$ is nonempty, and in fact contains
$D(c_1,c_2)$.  In particular, there is a path in $\overline{\bf T}_g$ from $\sigma_1$ to $\sigma_2$, lying completely in the boundary, which corresponds to
first pinching $c_2$, and then ``opening up" $c_1$.   Theorem \ref{T:corners} states that this path has a ``corner" at its intersection with $D(c_1,
c_2)$, and  is therefore not length minimizing.  In fact, the geodesic from
$\sigma_1$ to $\sigma_2$ intersects the boundary of $\overline{\bf T}_g$ only in its endpoints.
\end{Rem}

\begin{proof}[Proof of Theorem \ref{T:corners}]
The proof of this result may be modelled on that of Theorem \ref{T:interior} above.  We outline the approach here.  Assume, with the intent of
arriving at a contradiction, that $u:[0,1]\to\Delta_k$ is a geodesic with 
$
u(0)\in D(\ck)$,  $u(1)\in D(\ckprime)
$,
 $\ckprime\in \scckprime$, with $u(1)\not\in\overline D(\ck)$, and $k'\geq k$.  By changing endpoints, we may assume
$u$ is contained in 
$\overline D(\ck)\cup\overline D(\ckprime)
$.
It then follows from Theorem \ref{T:interior} that there is some $\hat x\in (0,1)$ such that
$$
u([0,\hat x))\subset D(\ck)\ ,\ u((\hat x,1])\subset D(\ckprime)\ ,\ u(\hat x)\in D({\bf c}_{(\ell)})\subset\overline D(\ck)\cap\overline D(\ckprime)\ .
$$
After renumbering, and by the assumption that $u(1)\not\in\overline D(\ck)$,  we may write
\begin{align*}
\ck&=\{c_1,\ldots, c_{k_1}, c_{k_1+1},\ldots, c_k\}\\ 
\ckprime&=\{c_1,\ldots, c_{k_1},
c_{k_1+1}^\prime,\ldots, c_{k'}^\prime\}\ ,
\end{align*}
for some $0\leq k_1<k$, and
 $c_i \neq c_j^\prime$ for any $i,j\geq k_1+1$.  Then by the exact same argument as was used in
Claim
\ref{C:dck} above, it follows that 
$$
{\bf c}_{(\ell)}=\{c_1,\ldots, c_{k_1}, c_{k_1+1},\ldots, c_k, c_{k_1+1}^\prime,\ldots, c_{k'}^\prime\}\ .
$$
Thus,  after reparametrization and forgetting about the first $k_1$ factors,  we have produced a geodesic
$$
u:[-1,1]\to\overline\HBbb^{k-k_1}\times\overline\HBbb^{k'-k_1}\times\Tbarprime
$$
satisfying:
\begin{align} 
u([-1,0))&\subset\left(\partial \HBbb\right)^{k-k_1}\times \HBbb^{k'-k_1}\times\Tbarprime\notag \\
u((0,1])&\subset\HBbb^{k-k_1}\times \left(\partial\HBbb\right)^{k'-k_1}\times\Tbarprime  \label{E:corner2}\\
u(0)&\in \left(\partial\HBbb\right)^{k-k_1}\times\left(\partial\HBbb\right)^{k'-k_1}\times\Tbarprime\notag
\end{align}
Using a rescaling argument as in the proof of Theorem \ref{T:interior}, we compare this geodesic to one in the model metric:
\begin{align*}
w:[-1,1]&\to\overline\HBbb^{k-k_1}\times\overline\HBbb^{k'-k_1}\\
w(-1)\in\left(\partial\HBbb\right)^{k-k_1}\times\HBbb^{k'-k_1} \ &,\ w(1)\in\HBbb^{k-k_1}\times\left(\partial\HBbb\right)^{k'-k_1}\ .
\end{align*}
Any such  geodesic with length bounded away from zero will keep $w(0)$ bounded away from $\left(\partial
\HBbb\right)^{k-k_1}\times\left(\partial\HBbb\right)^{k'-k_1}$.  As in the proof of Theorem
\ref{T:interior}, this leads to a contradiction to (\ref{E:corner2}).  
\end{proof}

We are now prepared to complete the 

\begin{proof}[Proof of Proposition \ref{P:wpgeodesic}] Since $\Tbar$ is a length space, $u$ is a geodesic. If $u$ has a point in $\T$, then by
Proposition
\ref{P:wolpert} (2),
$u$ is entirely contained in $\T$.  Suppose $u$ is entirely contained in $\partial\T$.  Then there is some $k\geq 1$ such that $u$ is
contained in $\Delta_k$ but not $\Delta_{k+1}$.  Therefore,  $u(x)\in D(\ck)$ for some $x$ and some $\ck\in\scck$.  By
equivariance, it follows that $u(x+1)=\gamma u(x)\in\dot\Delta_k$ as well.  By Theorem \ref{T:corners}, we must have that $\gamma
u(x)\in D(\ck)$.  But then $\gamma$ fixes $\ck$, contradicting the assumption of irreducibility.
\end{proof}


\section{Proof of Theorem \ref{T:equivariantgeodesic}}   \label{S:existence}


 Since $\gamma$ is pseudo-Anosov it follows that
$\twp(\gamma)>0$.  For if not we could find a sequence $\sigma_i\in \T$ such that $\dwp(\sigma_i,\gamma\sigma_i)\to 0$.  By
compactness we may assume that the corresponding sequence $[\sigma_i]\in\M$ converges to $[\sigma]$ in $\Mbar$.   If $[\sigma]\in\M$,
then by proper discontinuity of the action of
$\mcg$, $\gamma$ is periodic, which is a contradiction.  If $[\sigma]\in\partial\M$, then for $i$ large, $\gamma$ fixes the  collection of
short geodesics on $\sigma_i$.  This implies that $\gamma$ is reducible, which is again a
contradiction.

Note that it suffices to prove the result for some power $\gamma^k$ of the pseudo-Anosov:  If $A_{\gamma^k}$ is the unique
$\gamma^k$-equivariant geodesic, then since $\gamma A_{\gamma^k}$ is also a geodesic left invariant by $\gamma^k$, we must have $\gamma
A_{\gamma^k}=A_{\gamma^k}$, so that this is also a $\gamma$-equivariant geodesic.  With this understood, we may assume without loss of
generality that
$\gamma$ is in some finite index subgroup
$\mcg^\prime\subset\mcg$.  We will choose $\mcg^\prime$ to consist of {\bf pure} elements (cf.\ \cite[p.\ 3 and Corollary 1.8]{Iv}): these are
mapping classes that are either pseudo-Anosov or reducible.  Furthermore,  the reducible pure mapping classes have representatives that fix a
collection of simple closed curves pointwise, leave invariant the components of the complement, and on each component they are isotopic to
either the identity or a pseudo-Anosov.  This turns out to be a very useful construction, mainly for the following reason:

\begin{Lem} \label{L:pure}
Let $\mcg^\prime$ be a finite index subgroup of pure mapping classes.  Let $\Mprime$, $\Mbarprime$, and $[D(\ck)]'$ denote the
quotients of $\T$, $\Tbar$, and $D(\ck)$  by
$\mcg^\prime$.  Then $\Mbarprime$ is compact, $p':\T\to\Mprime$ is a covering space, and for each stratum $p':D(\ck)\subset\partial\T$,
the map
$D(\ck)\to [D(\ck)]'$ is also a covering.
\end{Lem}

\begin{proof}  Compactness of $\Mbarprime$ follows from the fact that $\Mbar$ is compact and $\mcg^\prime$ has finite index.  Moreover,
$\mcg^\prime$ acts properly discontinuously on $\T$ since $\mcg$ does.  The isotropy subgroups in $\mcg^\prime$ are trivial since this
$\mcg^\prime$ contains no periodic elements.  Hence, $\T\to\Mprime$ is a covering space.  The argument for the strata of $\partial\T$ is
similar.
\end{proof}

 Choose a length minimizing sequence $u_j:\RBbb\to
\T$ of smooth $\gamma$-equivariant paths with uniform modulus of continuity.  By convexity we have that $L(u_j)\to\twp(\gamma)$.  
 Let $v_j=p'\circ u_j$
denote the quotient map
$[0,1]\to\Mprime$.  By Ascoli's theorem and the compactness of $\Mbarprime$,  we may assume that after passing to a subsequence the $v_j$ converge uniformly
to a Lipschitz map
$v_\infty:[0,1]\to\Mbarprime$.  We again point out that there is no \`a priori reason for the existence of a convergent subsequence in $\Tbar$. For the sake
of clarity, we will start by considering the case where
$v_\infty$ does not lie entirely in the boundary.  Hence, there is a point $x_0$ such that $v_\infty(x_0)\in\Mprime$. Without loss of generality, we
reparametrize so that
$x_0=0$.

We now  modify the sequence by replacing each $u_j$ on $[0,1]$ by the \emph{geodesic} from $u_j(0)$ to $u_j(1)=\gamma u_j(0)$.
We then extend this modified sequence equivariantly to $\RBbb$.  We shall continue to denote it by $u_j$ and its projection to
$\Mprime$ by
$v_j$.  Note that since the replacement decreases length,
$u_j$ is still a length minimizing sequence and still has a uniform modulus of continuity.  We may assume, after again
passing to a subsequence, that $v_j$ converges uniformly to a Lipschitz map, again denoted by
$v_\infty:[0,1]\to\Mbarprime$ with
$v_\infty(0)=[\sigma_0]'\in\Mprime$.
By lower semicontinuity of lengths (cf.\ \cite[Proposition 1.20]{BH}),
$L(v_\infty)\leq \twp(\gamma)$.

By Proposition  \ref{P:wpgeodesic}, it then suffices to show that $v_\infty$ is the projection of a $\gamma$-equivariant path $\tilde v_\infty :
\RBbb\to\Tbar$.  To do this, set $S=v_\infty^{-1}(\partial\Mprime)$.  Since by Lemma \ref{L:pure}, $\T\to\Mprime$ is a covering map, if $S$ is
empty  the desired lift can be found.  So the closed set $S\subset (0,1)$ is the only possible obstruction to lifting, and we therefore assume
it is nonempty.  A posteriori, of course, the conclusion of Proposition  \ref{P:wpgeodesic} implies that this possibility does not, in fact, occur.

For simplicity, we will first also assume that the image of $S$ is contained in the top dimensional stratum of $\partial\Mprime$, i.e.
$v_\infty[0,1]\subset\Mprime\cup[\dot\Delta_1]'$,
 where $[\dot\Delta_1]'$ denotes the quotient of $\dot\Delta_1\subset\partial\T$ by $\mcg^\prime$. 

 Associated to the compact set $\displaystyle v_\infty(S)\subset[\dot\Delta_1]'$ we may choose $\kappa>0$ such that
\begin{enumerate}
\item[($\ast$)] The injectivity radius in $[\dot\Delta_1]'$ on each component of $v_\infty(S)$  is $>> \kappa$;
\item[($\ast\ast$)] 
 Let
$\widetilde {v_\infty(S)}$ denote the preimage of ${v_\infty(S)}$ in $\Tbar$.
 For any two nonisotopic simple closed curves $c, c'\in \scc$, the distance in $\Tbar$ between $\widetilde {v_\infty(S)}\cap D(c)$ and
$\widetilde {v_\infty(S)}\cap D(c')$  is
$>>
\kappa$.
\end{enumerate}
Both statements follow from the compactness of $v_\infty(S)$.

Recall that $[D(c)]'$ denotes the quotient of $D(c)\subset\partial\T$ by $\mcg^\prime$.  In particular,
$[\dot\Delta_1]'$ is the union of connected components  $[D(c)]'$ over isotopy classes of curves $c\in\scc$.
 Given $\delta>0$, we define a neighborhood $U_{c,\delta}$ of $[D(c)]'$ by the following condition:
$[\sigma]'\in U_{c,\delta}$ if $\ell_c(\sigma) <\delta$ for some representative $\sigma$ of $[\sigma]'$.  Also, let $\displaystyle
U_\delta=\cup\{ U_{c,\delta} : c\in\scc\}$. We will also need the following:

\begin{Lem}  \label{L:technical}
  Given any $\varepsilon>0$
 there is a $\delta>0$ such that if $[\sigma]'\in U_{c,\delta}$, then $[\sigma]'$ is within $\varepsilon$ of
$[D(c)]'$.  Moreover, there is a constant $\varepsilon_1>0$, depending only on  the choice of $\kappa$, with the following
significance: if $\varepsilon$ is sufficiently small compared to $\kappa$, 
if
$[\sigma]'\in U_{c,\delta}$ is within
$\varepsilon_1$ of
$v_\infty(S)\cap[D(c)]'$, and if
$\sigma_1,
\sigma_2$ are any two lifts of $[\sigma]'$ to $\Tbar$, then either $\dwp(\sigma_1,\sigma_2)\geq \kappa$, or $\dwp(\sigma_1,\sigma_2)\leq
2\varepsilon$.
\end{Lem}

\begin{proof}
The first statement follows from continuity of the length function with respect to the Weil-Petersson metric. 
Since $\mcg^\prime$ acts continuously  and consists of pure mapping classes, 
we may choose $\varepsilon_1$ sufficiently
small so that if $\sigma$ is within
$\varepsilon_1$ of  $\widetilde {v_\infty(S)}\cap D(c)$, and $\gamma\in\mcg^\prime$ satisfies $\gamma B_\kappa(\sigma)\cap B_\kappa(\sigma)\neq
\emptyset$, then $\gamma$
 acts trivially on $D(c)$.

With this choice of $\varepsilon_1$, let
$[\sigma]'$ satisfy the conditions of the lemma. Then by the first statement there is some $[\hat\sigma]'\in[D(c)]'$ such that
$\dwp([\sigma]',[\hat\sigma]')\leq
\varepsilon$.  Given any lift $\sigma_1$ of
$[\sigma]'$,  then since $\varepsilon<<\kappa$, ($\ast$) guarantees that 
$[\hat\sigma]'$ has a \emph{unique} lift
$\hat\sigma$ in a  $\kappa$-neighborhood of $\sigma_1$; moreover, $\dwp(\sigma_1, \hat\sigma)\leq \varepsilon$. As discussed above, if
$\sigma_2$ is another lift of $[\sigma]'$ in this neighborhood, i.e.\ $\sigma_2=\gamma\sigma_1$ for some $\gamma\in\mcg^\prime$, then
$\gamma$ acts trivially on $D(c)$.  In particular, $\gamma\hat\sigma=\hat\sigma$.  Hence,  
$\dwp(\sigma_2,\hat\sigma)\leq
\varepsilon$ as well.  The result follows.
\end{proof}

\setcounter{Claim}{0}

\begin{Claim}  \label{C:finite}
Let ${\mathcal O}=[0,1]\setminus S$.  Then
${\mathcal O}$ is a {\bf finite} union of relatively open intervals $I_\beta$.
\end{Claim}

\begin{proof}
If this were not the case, there would be a connected component $I_\beta=(x_\beta, y_\beta)$ of ${\mathcal O}$ with $v_\infty(I_\beta)$
arbitrarily short.  This follows from the uniform modulus of continuity.  We may assume that the length of  $v_\infty(I_\beta)$ is so small that
$v_\infty(I_\beta)$ is contained in a neighborhood of radius $\kappa/2$  of 
$v_\infty(S)\cap [D(c)]'$ for some $c$.  For
$j$ large, it follows that
$v_j(I_\beta)$ is contained in a $\kappa$-neighborhood of $v_\infty(S)\cap[D(c)]'$.  By ($\ast\ast$) it follows that there exist
simple closed curves
$c_j$, related to $c$ by elements of $\mcg^\prime$, such that $\ell_{c_j} u_j(x_\beta)$ and $\ell_{c_j} u_j(y_\beta)$ both tend to zero as
$j\to\infty$.  Since $u_j$ is a geodesic on $I_\beta$, it follows by Proposition \ref{P:wolpert} (3) that $\ell_{c_j} u_j(z)\to 0$ for any
$z\in I_\beta$.  But then
$v_\infty(I_\beta)\subset\partial \Mprime$, which is a contradiction.
\end{proof}

\begin{Claim} \label{C:neighborhood}
For  $\delta>0$ sufficiently small and  every  $\beta$, where $I_\beta=(x_\beta,y_\beta)$, there are curves
$c_{x_\beta}$ and $c_{y_\beta}$ and points $x_\beta\leq w_\beta< z_\beta\leq y_\beta$ such that:
\begin{enumerate}
\item   $v_\infty\left([x_\beta, w_\beta]\right)\subset \overline U_{c_{x_\beta},\delta}$
\item  $ v_\infty\left([z_\beta, y_\beta]\right)\subset\overline U_{c_{y_\beta},\delta} $
\item  $v_\infty\left([w_\beta,z_\beta]\right)\subset \Mprime\setminus U_\delta$
\par\smallskip\leftline{
If $I_\beta$ is of the form $[0,y_\beta)$ (resp.\ $(x_\beta,1]$), then we have}
\smallskip
\item   $v_\infty\left([z_\beta, y_\beta]\right)\subset\overline U_{c_{y_\beta},\delta}$
(resp.\ $v_\infty\left([x_\beta, w_\beta]\right)\subset \overline U_{c_{x_\beta},\delta}$)
\item  $v_\infty([0,z_\beta])\subset \Mprime\setminus U_\delta$
(resp.\ $v_\infty([w_\beta, 1])\subset \Mprime\setminus U_\delta$
\end{enumerate}
\end{Claim}

\begin{proof}
Since $v_\infty(I_\beta)\subset\Mprime$, we may find a $\delta$ satisfying (3) by compactness.  By choosing $\delta$ even smaller, and using the
uniform modulus of continuity of $v_\infty$, we may assume that $v_\infty(x_\beta,w_\beta)$ and $v_\infty(z_\beta,y_\beta)$ have very small
length compared to $\kappa$ and that $v_\infty(w_\beta)$ and $v_\infty(z_\beta)$ are on the boundary of the ball $U_{c_{x_\beta},\delta}$.  Thus, by the
same argument as in the proof of Claim
\ref{C:finite}, there exist
$c_j$, related to some
$c_{x_\beta}$ by elements of $\mcg^\prime$, such that
$
\ell_{c_j} u_j(x_\beta)\to 0$, and $\ell_{c_j} u_j(w_\beta)\to \delta
$
as $j\to \infty$.  By convexity, it follows that for any point $z\in (x_\beta, w_\beta)$, $\displaystyle\limsup_{j\to\infty}\ell_{c_j}
u_j(z)\leq\delta$.  It follows that $v_\infty(z)\in \overline U_{c_{x_\beta},\delta}$.  The other parts follow similarly.
\end{proof}

\begin{Claim} \label{C:betalift}
For $\varepsilon_2>0$ sufficiently small and $j$ sufficiently large there is a lift $\tilde v_{\infty,\beta}$ of $v_\infty$ on each $I_\beta$
such that $\sup_{x\in I_\beta}\dwp (\tilde v_{\infty,\beta}(x), u_j(x))\leq 2\varepsilon_2$.  Moreover, if $\delta>0$ and $x_\beta, y_\beta$ are
the endpoints of $I_\beta$, we can arrange that $u_j(x_\beta)$ and $u_j(y_\beta)$ are in $U_\delta$.
\end{Claim}

\begin{proof}
Fix $\varepsilon_2$ small relative to $\kappa$ and $\varepsilon_1$.  
Choose $\delta$ small relative to $\varepsilon_2$ so that the conclusion of Lemma \ref{L:technical} (with $\varepsilon_2$ as $\varepsilon$)
holds.  We may choose
$\delta$ possibly smaller so that the conclusion of Claim
\ref{C:neighborhood} holds.  On the complement of
$U_\delta$ there is a uniform lower bound on the injectivity radius of $\Mprime$.  
Choose $\varepsilon\leq \varepsilon_2$  small compared to this.  For $j$ sufficiently large $v_j$ is within $\varepsilon$ of $v_\infty$
uniformly on 
$[0,1]$.   Thus, for such a $j$, $v_\infty$ lifts on $I_\beta$ to a path
$\tilde v_{\infty,\beta}$ that is within
$\varepsilon$ of
$u_j$  at some point, and hence every point in $(w_\beta, z_\beta)$.  Now by Lemma \ref{L:technical}, any two lifts of
$v_\infty$ on $(x_\beta, w_\beta)$ are either uniformly within $2\varepsilon_2$ of each other, or they are separated by a distance at least
$\kappa$. Since we assume $\varepsilon\leq\varepsilon_2 << \kappa$, we see that the distance between the lift $\tilde v_{\infty,\beta}$ and $u_j$
on
$(x_\beta, w_\beta)$ is at most $2\varepsilon_2$.  The same reasoning applies to $(z_\beta, y_\beta)$.  Finally, since by Claim \ref{C:finite}
there are only finitely many intervals
$I_\beta$, we may choose $j$ uniformly.  In particular, we can guarantee the second statement.
\end{proof}

Recall that ${\mathcal O}$ is the finite union of intervals $ I_\beta$.  For any $\varepsilon_2>0$ sufficiently small we have defined a lift $\tilde
v_{\infty,{\mathcal O}}:{\mathcal O}\to\T$ of
$v_\infty$ such that 
$ \sup_{x\in I_\beta}\dwp (\tilde v_{\infty,\beta}(x), u_j(x))\leq 2\varepsilon_2$, where $j$ is chosen large according to
$\varepsilon_2$.  We first point out  that:

\begin{Claim} \label{C:equivariant}
 $\tilde v_{\infty,{\mathcal O}}(1)=\gamma\tilde v_{\infty,{\mathcal O}}(0)$.
\end{Claim}

\begin{proof}
The two points are both lifts of $[\sigma_0]'\in\Mprime$, which lies in the complement of $U_\delta$.  By the equivariance of $u_j$ and the
construction in Claim \ref{C:betalift}, they are also within
$2\varepsilon$ of each other.  Since $\varepsilon$ has been chosen much smaller
than the injectivity radius of
$\Mprime$ on the complement of
$U_\delta$, the points must agree.
\end{proof}

 By the uniform
modulus of continuity, $\tilde v_{\infty, {\mathcal O}}$ extends to a map $\tilde v_{\infty,\overline {\mathcal O}}:\overline {\mathcal O}\to\Tbar$ which is still a lift of
$v_\infty$.

\begin{Claim} \label{C:continuity}
 $\tilde v_{\infty,\overline {\mathcal O}}$ is  continuous.
\end{Claim}

\begin{proof}
Consider intervals $I_\beta=(x_\beta,y_\beta)$ (or $[0,y_\beta)$) and $I_{\beta'}=(x_{\beta'},y_{\beta'})$  (or $(x_{\beta'}, 1]$) with
$y_\beta=x_{\beta'}$.  Then
$\tilde v_{\infty,\overline {\mathcal O}}(y_\beta)$ and $\tilde
v_{\infty,\overline {\mathcal O}}(x_{\beta'})$ are two lifts of $v_\infty(y_\beta)$.  Moreover, since both points are within $2\varepsilon_2$ of
$u_j(y_\beta)$, it follows that they are within $4\varepsilon_2$ of each other.  On the other hand, they are lifts of points in $v_\infty(S)$, and
since
$\varepsilon_2$ is small compared to $\kappa$, they must coincide by ($\ast$).
\end{proof}

Now for each interval $J=(a,b)\subset [0,1]\setminus\overline {\mathcal O}$, we have $v_\infty(J)\subset[D(c)]'$ for some $c$.  Since $D(c)\to[D(c)]'$ is a
covering, there is a unique lift  $\tilde v_{\infty, J}$ of $v_\infty $ on $J$ such that $\lim_{x\downarrow a}\tilde
v_{\infty,J}(x)=\tilde v_{\infty,\overline {\mathcal O}}(a)$.

\begin{Claim}  \label{C:something}
 $\sup_{x\in J}\dwp(\tilde v_{\infty, J}(x), u_j(x))\leq 2\varepsilon_2$.
\end{Claim}

\begin{proof}
By ($\ast\ast$), the fact that $v_j$ is uniformly within $\varepsilon<<\kappa$ of $v_\infty$, and the second statement in Claim \ref{C:betalift},
it follows that there is  some representative
$c_j$ of $c$ such that 
$
\ell_{c_j} u_j(a)$ and $  \ell_{c_j} u_j(b)
$
are both less than $\delta$.  By convexity, this holds on the entire interval $J$.  Hence, $u_j$ is within $\varepsilon_2$ of $\widetilde
{v_\infty(S)}\cap D(c_j)$.   As in the proof of  Lemma \ref{L:technical}, a lift $\tilde v_{\infty, J}$ is either uniformly within $\varepsilon_2$
of
$u_j$ or uniformly a distance at least $\kappa$ away.  Since the latter case does not hold at the initial point $a$, the result follows.
\end{proof}

Given Claim \ref{C:something}, it follows exactly as in Claim \ref{C:continuity} that $\lim_{x\uparrow b}\tilde
v_{\infty,I}(x)=\tilde v_{\infty,\overline {\mathcal O}}(b)$.  Continuing this construction for each connected component $J$ of
$[0,1]\setminus\overline{\mathcal O}$, we finally end up with a lift $\tilde v_\infty$ of
$v_\infty$ on $[0,1]$.  By Claim \ref{C:equivariant}, we can extend $\tilde v_\infty$ to a continuous equivariant path $\tilde v_\infty :
\RBbb\to\Tbar$.  By Claim \ref{C:finite} it follows that  $L(\tilde v_\infty)=L(v_\infty)\leq\twp(\gamma)$.  By Proposition
\ref{P:wpgeodesic}, $\tilde v_\infty$ is a $\gamma$-equivariant geodesic in $\T$, and this completes the existence part of the proof of Theorem
\ref{T:equivariantgeodesic} in the case where $v_\infty([0,1])\subset\Mprime\cup[\dot\Delta_1]'$.

The general case will follow from an inductive argument.   First, note that the case where 
$v_\infty([0,1])\subset[\dot\Delta_s]'\cup[\dot\Delta_{s+1}]'$, $s\geq 1$, follows exactly as above.  We also point out that in the
argument given above we may assume without loss of generality that $v_\infty(0)$  lies in $[\dot\Delta_s]'$.  This is because the
argument may be easily adapted to a minimizing sequence of once-broken geodesics rather than geodesics. Suppose now that
$$
v_\infty([0,1])\subset[\dot\Delta_s]'\cup[\dot\Delta_{s+1}]'\cup\cdots\cup[\dot\Delta_{k}]'\ ,
$$
where $0\leq s\leq k$.  We induct on $k-s$.  More precisely, we prove the following statement:

\medskip
\noindent
($\ast\ast\ast$)\ \it Given a path $v_\infty=\lim v_j$ with $v_j=p'\circ u_j:[0,1]\to\Mbarprime$, where $u_j:[0,1]\to \Tbar$ are geodesics with
uniform modulus of continuity such that
$$
v_\infty([0,1])\subset[\dot\Delta_s]'\cup[\dot\Delta_{s+1}]'\cup\cdots\cup[\dot\Delta_{k}]'\ ,\ v_\infty(0) \in[\dot\Delta_s]'\ ,
$$
then for $\varepsilon>0$ sufficiently small there is a lift $\tilde v_\infty$ of $v_\infty$ such that 
$$
\max_{x\in[0,1]}\dwp (u_j(x),\tilde v_\infty(x)) <\varepsilon\ .
$$
Furthermore, if the sequence satisfies $u_j(1)=\gamma u_j(0)$ for some $\gamma\in\mcg$, then the $\tilde v_\infty$ can be chosen so that $\tilde
v_\infty(1)=\gamma\tilde v_\infty(0)$.\rm
\medskip

The statement ($\ast\ast\ast$) has already be proven in the case $k-s=1$.  Inductively assume the statement to be true for all integers $< k-s$, and let
$v_\infty$ be as in ($\ast\ast\ast$).  Let $S=v_\infty^{-1}[\dot\Delta_k]'$.  If $S=\emptyset$ the statement holds by the inductive
hypothesis.  Thus we may assume that $S$ is a nonempty closed subset of $(0,1)$.  Next, associated to the compact set
$v_\infty(S)\subset[\dot\Delta_k]'$ we can choose $\kappa>0$ such that:
\begin{enumerate}
\item[($\ast_k$)] The injectivity radius in $[\dot\Delta_k]'$ on each component of $v_\infty(S)$  is $>> \kappa$;
\item[($\ast\ast_k$)] 
 Let
$\widetilde {v_\infty(S)}$ denote the preimage of ${v_\infty(S)}$ in $\Tbar$.
 For any two nonisotopic  $\ck, \ck'\in \scck$, the distance in $\Tbar$ between $\widetilde {v_\infty(S)}\cap D(\ck)$ and
$\widetilde {v_\infty(S)}\cap D(\ck')$  is
$>>
\kappa$.
\end{enumerate}
As before, both statements can be satisfied by the compactness of $v_\infty(S)$.

Given a collection $\ck=\{c_1,\ldots,c_k\}\in\scck$  we define a neighborhood $U_{\ck,\delta}\subset[\Delta_s]'$ of $[D(\ck)]'$ by the following
condition:
$[\sigma]'\in U_{\ck,\delta}$ if $\ell_{c_i}(\sigma) <\delta$, $i=1,\ldots, k$, for some representative $\sigma$ of $[\sigma]'$.  Also, let
$\displaystyle U_\delta=\cup\{ U_{\ck,\delta} : \ck\in\scck\}$.
Then a version of Lemma \ref{L:technical} holds in this situation as well, \emph{mutatis mutandi}.  Also, if we set ${\mathcal O}=[0,1]\setminus S$,
then $\mathcal O$ is a finite union of intervals as in Claim \ref{C:finite}.  The proof of Claim \ref{C:neighborhood} also goes through with no change. 
For convenience, we state the appropriate version:
\setcounter{Claim}{1}
\begin{Claim}[$k$] \label{C:neighborhoodk}
For  $\delta>0$ sufficiently small and  every  $\beta$, where $I_\beta=(x_\beta,y_\beta)$, there are curves
${\ck}_{x_\beta}$ and ${\ck}_{y_\beta}$ and points $x_\beta\leq w_\beta< z_\beta\leq y_\beta$ such that:
\begin{enumerate}
\item   $v_\infty\left([x_\beta, w_\beta]\right)\subset \overline U_{{\ck}_{x_\beta},\delta}$
\item  $ v_\infty\left([z_\beta, y_\beta]\right)\subset\overline U_{{\ck}_{y_\beta},\delta} $
\item  $v_\infty\left([w_\beta,z_\beta]\right)\subset [\dot\Delta_s]'\cup[\dot\Delta_{s+1}]'\cup\cdots\cup[\dot\Delta_{k-1}]'\setminus U_\delta$
\par\smallskip\leftline{
If $I_\beta$ is of the form $[0,y_\beta)$ (resp.\ $(x_\beta,1]$), then we have}
\smallskip
\item   $v_\infty\left([z_\beta, y_\beta]\right)\subset\overline U_{{\ck}_{y_\beta},\delta}$
(resp.\ $v_\infty\left([x_\beta, w_\beta]\right)\subset \overline U_{{\ck}_{x_\beta},\delta}$)
\item  $v_\infty([0,z_\beta])\subset \Mprime\setminus U_\delta$
(resp.\ $v_\infty([w_\beta, 1])\subset [\dot\Delta_s]'\cup[\dot\Delta_{s+1}]'\cup\cdots\cup[\dot\Delta_{k-1}]'\setminus U_\delta$
\end{enumerate}
\end{Claim}

\begin{Claim}[$k$] \label{C:betaliftk}
For $\varepsilon_2>0$ sufficiently small and $j$ sufficiently large there is a lift $\tilde v_{\infty,\beta}$ of $v_\infty$ on each $I_\beta$
such that $\sup_{x\in I_\beta}\dwp (\tilde v_{\infty,\beta}(x), u_j(x))\leq 2\varepsilon_2$.  Moreover, if $\delta>0$ and $x_\beta, y_\beta$ are
the endpoints of $I_\beta$, we can arrange that $u_j(x_\beta)$ and $u_j(y_\beta)$ are in $U_\delta$.
\end{Claim}

\begin{proof}
Again choose  $\varepsilon_2$ and $\delta$ as in Claim \ref{C:betalift}.  Then choose $\varepsilon\leq \varepsilon_2$ less than the injectivity radius
of the strata $\dot\Delta_{k'}\setminus U_\delta$, $s\leq k'<k$.  By the inductive hypothesis, if $j$ is sufficiently large there is a lift $\tilde
v_{\infty,\beta}$ of $v_\infty$ on $I_\beta$ that is within $\varepsilon$ of $u_j$ at every point in $(w_\beta, z_\beta)$.  Now by the analogue of Lemma
\ref{L:technical},  any two lifts of
$v_\infty$ on $(x_\beta, w_\beta)$  or $(z_\beta, y_\beta)$ are either uniformly within $2\varepsilon_2$ of each other, or they are separated by a
distance at least
$\kappa$. Since we assume $\varepsilon\leq\varepsilon_2 << \kappa$, we see that the distance between the lift $\tilde v_{\infty,\beta}$ and $u_j$
on
$(x_\beta, w_\beta)$ is at most $2\varepsilon_2$.   As before,  since 
there are only finitely many intervals
$I_\beta$, we may choose $j$ uniformly, and the second statement follows.
\end{proof}

In this way we have constructed a lift $\tilde v_{\infty,{\mathcal O}}:{\mathcal O}\to\Tbar$.  Note that because of the choice of $\varepsilon$ in Claim
\ref{C:betaliftk} $(k)$,  the exact same statement Claim \ref{C:equivariant} applies to $\tilde v_{\infty,{\mathcal O}}$ as well.
Now the remaining parts of the argument, namely, Claims \ref{C:continuity} and \ref{C:something}, follow exactly as before to give the desired lift
$\tilde v_\infty$.  This proves  ($\ast\ast\ast$) for $k-s$, which by induction then holds in general.
By Proposition
\ref{P:wpgeodesic}, $\tilde v_\infty$ is a $\gamma$-equivariant geodesic in $\T$, and this completes the existence part of the proof of Theorem
\ref{T:equivariantgeodesic}.    Uniqueness follows from the geodesic convexity  and negative curvature of $\T$.


\section{Classification of Weil-Petersson Isometries}  \label{S:classification}


The existence of  equivariant geodesics for pseudo-Anosov mapping classes allows for the precise 
classification of  Weil-Petersson isometries in terms of translation length (\ref{E:translationlength}) that we have given in Table 1.  First, let us
clarify the terminology used there:  
$\gamma\in\mcg$ is {\bf pseudoperiodic} if it is either periodic, or it is reducible and periodic on the reduced components;  it is
{\bf strictly pseudoperiodic} if it is pseudoperiodic but
not periodic.  Furthermore, we say that $\gamma$ is  {\bf semisimple}  if there is  $\sigma\in\T$ such that
$\twp(\gamma)=\dwp(\sigma,\gamma\sigma)$.

\begin{Thm} \label{T:classification}
An infinite order mapping class  is semisimple if and only if it is irreducible.
More precisely, any $\gamma\in\mcg$
belongs to exactly one of the four classes characterized in Table 1.
\end{Thm}

\smallskip
\noindent    A slightly different picture emerges if one considers the action of $\mcg$ on the completion $\Tbar$ of $\T$.  There, every $\gamma\in\mcg$
is semisimple, although  the action of $\mcg$ on $\Tbar$ is, of course, no longer properly discontinuous (indeed, according to \cite{KL} it cannot be).
It is worth mentioning that while we have defined pseudo-Anosov's in the usual way in terms of measured foliations, the only property that we use
in the proofs of Theorems \ref{T:equivariantgeodesic} and \ref{T:classification} is that these mapping classes have infinite order and are
irreducible.  In particular, the description given in Table 1 is independent of Thurston's classification provided we replace the entry
``pseudo-Anosov" by ``infinite irreducible."

To begin the proof, note that 
the following is a consequence of the proper discontinuity of the action of the mapping class group that we have already implicitly used 
in the proof of Lemma \ref{L:technical}:

\begin{Lem}  \label{L:discontinuous}
Let $\sigma\in\Tbar$.  Then there is a neighborhood $U$ of $\sigma$ such that if $\gamma\in\mcg$ is such that $\gamma U\cap
U\neq\emptyset$, then $\gamma$ is pseudoperiodic.
\end{Lem}

\begin{Lem} \label{L:pseudoperiodic}
Let $\gamma\in\mcg$.  Then $\twp(\gamma)=0$ if and only if $\gamma$ is pseudoperiodic.
\end{Lem}

\begin{proof}
Clearly, we need only show that $\twp(\gamma)=0$ implies pseudoperiodic, the other direction following from the fact that a pseudoperiodic
element has a fixed point in $\Tbar$.  Let
$\sigma_j\in\T$ be a sequence such that $\dwp(\sigma_j,\gamma\sigma_j)\to 0$.  We may assume that $[\sigma_j]$ converges to some
$[\sigma]\in\Mbar$.  Hence, we can find $h_j\in\mcg$ such that $h_j\sigma_j$ converges to some lift $\sigma\in\Tbar$ of
$[\sigma]$.  Let $U$ be a neighborhood of $\sigma$ as in Lemma \ref{L:discontinuous}.  For $j$ large, $h_j\sigma_j$ and $h_j\gamma\sigma_j$
are in $U$.  Then by the conclusion of the lemma, 
$h_j\gamma h_j^{-1}$ is pseudoperiodic; hence, so is $\gamma$.
\end{proof}

 \begin{Lem}  \label{L:pseudo}
If there is a complete, nonconstant Weil-Petersson geodesic in $\T$ that is equivariant with respect to a mapping class $\gamma\in \mcg$, then 
$\gamma$ is pseudo-Anosov.
\end{Lem}
\begin{proof}
By geodesic convexity and negative curvature, $\gamma$ cannot be periodic.  We claim that $\gamma$ cannot be reducible either. To prove the claim,
let
$c\in\scc$.   By Wolpert's result Proposition \ref{P:wolpert} (3), if $u:\RBbb\to \T$ is a Weil-Petersson
geodesic, then
$f(x)=\ell_c(u(x))$ is a strictly convex function.  On the other hand, if
$\gamma(c)=c$ up to isotopy, and if $u$ is equivariant with respect to $\gamma$, then $f(x)$ would also be periodic, which is absurd.
\end{proof}

\begin{proof}[Proof of Theorem \ref{T:classification}]
The first row of Table 1 follows from Lemma \ref{L:pseudoperiodic} and the fact that $\gamma$ has a fixed point in $\T$ if and only if
$\gamma$ is periodic.  If $\gamma$ is pseudo-Anosov, then as a consequence of Theorem \ref{T:equivariantgeodesic}, $\twp(\gamma)$ is
attained along an equivariant geodesic, so pseudo-Anosov's are semisimple.
Conversely, suppose the infimum is attained at $\sigma\in\T$, but $\gamma\sigma\neq\sigma$.  
Then we argue as in Bers \cite{B} (see also, \cite[p.\ 81]{BGS}): by
Proposition \ref{P:wolpert} (2) there is a  Weil-Petersson geodesic $u$ from $\sigma$ and $\gamma\sigma$.  Similarly, there is a geodesic $\beta$
from
$\sigma$ to $\gamma^2\sigma$.  We claim that $\beta=u\cup\gamma u$. For if not, and if $\sigma'$
denotes the point on $u$ half the distance from $\sigma$ to $\gamma\sigma$, then the distance from $\sigma'$ to $\gamma\sigma'$ would be
strictly less than
$\dwp(\sigma,\gamma\sigma)=\twp(\gamma)$, which is a contradiction.  In this way, we see that $u$ may be extended to a complete $\gamma$-equivariant
geodesic.   Lemma
\ref{L:pseudo} then implies that
$\gamma$ is pseudo-Anosov.  A process of elimination fills in the remaining item of the table.
\end{proof}


\section{Proof of Theorem \ref{T:proper}}  \label{S:proper}


We begin this section by studying the asymptotic behavior of Weil-Petersson geodesics.  The guiding idea is 
that the ideal boundary of $\T$ with respect to the Weil-Petersson metric should not be too different from the
Thurston boundary of projective measured foliations.  While we shall not attempt a full description here, the result we obtain in this direction
nevertheless indicates a significant link between  the Weil-Petersson geometry of the action of  $\mcg$
on
$\T$ on the one hand, and the purely topological action of $\mcg$ on $\PMF$ on the other.

\begin{Def} Let $\alpha, \alpha' :\RBbb\to\T$ be paths parametrized by arc length.  We say that $\alpha$ and $\alpha'$ {\bf diverge} if the
function
$
(t,s)\mapsto \dwp(\alpha(t),\alpha'(s))
$
is proper.
\end{Def}
 
   Our goal is to prove the following

\begin{Thm} \label{T:diverge}
Let $A_{\gamma}$ and  $A_{\gamma'}$ be the axes for independent
pseudo-Anosov mapping classes $\gamma$ and $ \gamma'$.  Then
$A_{\gamma}$ and $A_{\gamma'}$ diverge.
\end{Thm}

\noindent  As a first step, we have

\begin{Lem}  \label{L:thurstonlimit}
Let $A_\gamma$ be the axis of a pseudo-Anosov mapping class $\gamma$, and  let
$\{F_+,F_-\}$ denote the stable and unstable  foliations of $\gamma$.  If $A_\gamma(t_j)\to [F]$ in the Thurston topology for some sequence
$t_j\to\infty$, then $[F]\in\left\{[F_+],[F_-]\right\}$.
\end{Lem}

\begin{proof}
By definition, there is a fixed length $L=\twp(\gamma)$ such that $\dwp(A_\gamma(t_j),\gamma A_\gamma(t_j))=L$ for all $j$.   By
equivariance and the parametrization by arc length,
$\gamma A_\gamma(t_j)=A_\gamma(t_j+L)$. Since the
projection of
$A_\gamma$ to
$\M$ is a closed curve, it lives in a fixed compact subset.  Hence, along the path
$A_\gamma(t)$ the Weil-Petersson and Teichm\"uller metrics are uniformly quasi-isometric.  Thus, for some constant $K$,
$\dt(A_\gamma(t_j),\gamma A_\gamma(t_j))\leq K$ for all $j$.  For any simple closed curve $c$ on $\Sigma_{g,n}$, it follows from \cite[Lemma
3.1]{W2} that 
\begin{equation} \label{E:quasiconformal} 
e^{-2K}\ell_c (A_\gamma(t_j))\leq
\ell_c(\gamma A_\gamma(t_j))\leq e^{2K}\ell_c(A_\gamma(t_j))
\end{equation} 
for all $j$.  The result now follows easily: recall from Section \ref{S:mcg} the notion of a minimal foliation and the associated definitions
(\ref{E:top}) and (\ref{E:gr}).
If $[F]\not\in\Fmin$, then it is easily deduced from  (\ref{E:quasiconformal}) that $i(F,c)=0 \iff i(\gamma
F,c)=0$.  Hence, $[\gamma F]\in\ZF$, contradicting (\ref{E:zf}).  On the other hand, if
$[F]\in\Fmin$, then choose a sequence
$c_j\in\scc$  such that
$\ell_{c_j}(A_\gamma(t_j))$ is uniformly bounded for all $j$.  
We may assume there are numbers $r_j\leq 1$ and  $G\in\MF$ so that $r_jc_j\to G$ in the Thurston topology.
It then follows as in \cite[Lemma 2.1]{DKW} that $i(F,G)=0$. By
(\ref{E:quasiconformal}), we also have that
$\ell_{c_j}(\gamma A_\gamma(t_j))$ is uniformly bounded, so that by the same argument $i(\gamma F, G)=0$.  Since we assume $F$ is
minimal, this implies $[\gamma F]\in \WF$.  Hence, by (\ref{E:wf}), $[F]\in\left\{[F_+],[F_-]\right\}$.
\end{proof}

\begin{Lem}  \label{L:minimal}
Let $\sigma_j$, $\sigma_j^\prime$ converge in the Thurston topology to projective measured foliations $[F]$ and $[F']$, respectively.  Assume that $F$ is
minimal and that $[\sigma_j]$ and $[\sigma_j^\prime]$ lie in a fixed compact subset of $\M$.  Also, we suppose that there is a fixed $D\geq
0$ such that $\dwp(\sigma_j,\sigma_j^\prime)\leq D$ for all $j$.  Then $[F']\in\WF$.
\end{Lem}

\begin{proof}
Let $u_j$ denote the Weil-Petersson geodesic in $\T$ from $\sigma_j$ to $\sigma_j^\prime$, and let $v_j$ denote its projection to $\M$.  As in Section
\ref{S:existence}, we may assume (after passing to a subsequence) that the $v_j$ converge to a piecewise geodesic $v_\infty:[0,1]\to\Mbar$.  More
precisely, there are values
$0=x_0<x_1\leq x_2\leq \cdots\leq x_{N-1}<x_N=1$ with the following properties:  (1) $v_\infty\bigr|_{(x_i,x_{i+1})}$ is either a geodesic in $\M$ or
is entirely contained in $\partial\M$; and (2) if $v_\infty\bigr|_{[x_i,x_{i+1}]}\subset\partial\M$ then for $j$ sufficiently large there is a
sequence
$c_j^{(i)}$ of simple closed curves on $\Sigma_{g,n}$ whose lengths with respect to both $u_j(x_i)$ and $u_j(x_{i+1})$ converge to zero as
$j\to\infty$.  We now proceed inductively:  For any $x\in[0,x_1)$ and any sequence such that $u_j(x)\to [F_x]$ in the Thurston topology, we
claim that $F_x$ is minimal and $i(F,F_x)=0$.  For since $v_\infty\bigr|_{[0,x]}$ is contained in a compact subset of $\M$, on which  the
Weil-Petersson and Teichm\"uller metrics are therefore uniformly quasi-isometric, the claim follows exactly as in the proof of Lemma
\ref{L:thurstonlimit} above.  Now choose
$x<x_1$ and $x_2<y$ sufficiently close so that the lengths of $c^{(1)}_j$ with respect to both $u_j(x)$ and $u_j(y)$ are bounded.  If we choose
$r_j\leq 1$ such that $r_jc_j\to G$, then as in the proof of Lemma \ref{L:thurstonlimit} above, $i(F,G)=0$.  Similarly, if $u_j(y)\to [F_y]$,
then
$i(G,F_y)=0$, and so $F_y\in\WF$.  Continuing in this way, we see that $F'\in\WF$.
\end{proof}

\begin{proof}[Proof of Theorem \ref{T:diverge}]
Suppose $A_{\gamma}$ and $A_{\gamma'}$ do not diverge.  Then there is a constant $D\geq 0$ and an unbounded sequence $(t_j, s_j)$ such that
$\dwp(A_\gamma(t_j), A_{\gamma'}(s_j))\leq D$.  Necessarily, $t_j$ and $s_j$ are both unbounded, so we may assume that $t_j, s_j\to\infty$. Let
$\{F_+,F_-\}$ and $\{F_+^\prime,F_-^\prime\}$ denote the stable and unstable foliations of $\gamma$ and $\gamma'$, respectively. Since $\gamma$
and $\gamma'$ are independent,  all four of these foliations are mutually topologically distinct  (cf.\
\cite[Corollary 2.6]{McP}).
 After perhaps passing to  a subsequence, we may assume that there are measured foliations $F$ and $F'$ such that
$A_\gamma(t_j)\to [F]$ and
$A_{\gamma'}(s_j)\to [F']$ in the Thurston topology.  By Lemma \ref{L:thurstonlimit}, $[F]\in\left\{[F_+],[F_-]\right\}$ and
$[F']\in\{[F_+^\prime],[F_-^\prime]\}$.
 Let $\sigma_j=A_\gamma(t_j)$ and $\sigma_j^\prime=A_{\gamma'}(s_j)$.  Then since $\dwp(\sigma_j, \sigma_j^\prime)\leq D$, and  $[\sigma_j],
[\sigma_j^\prime]$ lie on the closed geodesics $[A_\gamma]$ and $[A_{\gamma'}]$, and the stable and unstable foliations of a pseudo-Anosov are
minimal,   the hypotheses of Lemma
\ref{L:minimal} are satisfied.  Thus, $F$ and $F'$ are topologically equivalent, which is a contradiction.  This completes the proof.
\end{proof}

  Let $H\subset\mcg$ be a
finitely generated subgroup.  To each set of generators $\mathcal G$ we associate a function on $\T$:
\begin{equation}  
\distortion(\sigma) = \max\left\{ \dwp(\sigma,\gamma\sigma) : \gamma\in{\mathcal G}\right\}\ .
\end{equation}

\begin{Def}[cf.\ {\cite[\S 2]{KS2}}] \label{D:proper}  
A finitely generated subgroup $H\subset \mcg$ is {\bf proper} (or acts {\bf properly} on $\T$) if there exists a generating set $\mathcal G$ of $H$ with
the property that for every
$M\geq 0$, the set $\{\sigma\in\T : \distortion(\sigma)\leq M\}$ is bounded.
\end{Def}

\noindent   Clearly, the condition in the definition of a proper subgroup is independent of the choice of generating set.   For complete manifolds
of nonpositive curvature, the existence of two  hyperbolic isometries with nonasymptotic axes is sufficient to prove that a subgroup of isometries
is proper.  Theorem \ref{T:proper} shows  that this works for mapping class groups as well, where hyperbolic is replaced with pseudo-Anosov.
  The approach
taken below to prove this result uses the relationship, demonstrated above, between the axes of  pseudo-Anosov's and their corresponding
stable and unstable foliations.  We first need one more

\begin{Lem} \label{L:temp}
Let $A_\gamma$ be the axis of a pseudo-Anosov $\gamma$ and $\sigma_i\in\T$ any sequence.  Then as $i\to\infty$,
$$
 \distwp(\sigma_i,A_\gamma)\to +\infty\quad \Longrightarrow\quad
\dwp(\sigma_i, \gamma\sigma_i)\to +\infty\ .
$$
\end{Lem}

\begin{proof}
Choose $\sigma\not\in A_\gamma$, and let $\sigma_0\in A_\gamma$ such that $\distwp(\sigma, A_\gamma)=\dwp(\sigma,\sigma_0)$.  Let $\alpha$ be the
geodesic from $\sigma_0$ to $\sigma$, parametrized by arc length.  Now $f(t)=\dwp(\alpha(t),\gamma\alpha(t))$ is strictly convex with $f'(0)=0$,
since $f(0)=\twp(\gamma)$ is the minimum of $f$.  Moreover, since
$A_\gamma$ projects to a closed curve in $\M$, there is a neighborhood of
$A_\gamma$ on which the curvature is bounded above by a strictly negative constant.  Therefore, for  $t_\gamma>0$ sufficiently small depending
only $A_\gamma$, there is
$\varepsilon_\gamma>0$, depending only on $t_\gamma$ and $A_\gamma$, such that $f'(t_\gamma)\geq \varepsilon_\gamma$.  So if
$L=\dwp(\sigma,\sigma_0)$, then
\begin{equation}  \label{E:distance}
\dwp(\sigma,\gamma\sigma)=f(L)\geq \varepsilon_\gamma(L-t_\gamma)=\varepsilon_\gamma(\dwp(\sigma,\sigma_0)-t_\gamma)\ .
\end{equation}
  Applying this to the sequence $\sigma_i$ completes the proof.
\end{proof}

\begin{proof}[Proof of Theorem \ref{T:proper}]
Let $H\subset\mcg$ be sufficiently large.  Then by definition $H$ contains two independent pseudo-Anosov's $\gamma$ and $\gamma'$.
Since the condition  of being proper does not depend on the choice of generating set, we may include the elements $\gamma, \gamma'$ in $\mathcal
G$.  If
$H$ is not proper then there is a number $M\geq 0$ and  an unbounded sequence $\sigma_i\in\T$ such that $\distortion(\sigma_i)\leq M$.  In
particular,
$\dwp(\sigma_i,\gamma\sigma_i)$ and
$\dwp(\sigma_i,\gamma'\sigma_i)$ are both bounded by $M$ for all $i$.
By Lemma \ref{L:temp}, it follows that both $\distwp(\sigma_i, A_\gamma)$ and $\distwp(\sigma_i, A_{\gamma'})$ are bounded.  Since $\sigma_i$ is
unbounded in $\T$, this implies that $\dwp(A_\gamma(t_i), A_{\gamma'}(s_i))$ is
bounded along some unbounded sequence $(t_i, s_i)$.  But this contradicts Theorem \ref{T:diverge}.
\end{proof}


\section{Further Results}       \label{S:further}


This section is independent of the rest of the paper.
In it we discuss the convergence of the heat equation as a method for finding equivariant geodesics.  The result, Theorem
\ref{T:convergence}, relies on the existence Theorem \ref{T:equivariantgeodesic}.  It would be useful to have a direct proof of convergence that
would circumvent the detailed discussion in Section \ref{S:existence}, though at present it is not clear how to do this.  At the end of this section we
also point out additional consequences of the existence of equivariant geodesics.

Let $\gamma\in\mcg$ and $u:\RBbb\to\T$, equivariant with respect to $\gamma$.  The length of $u$ is, by our convention, 
\begin{equation} \label{E:length}
L(u)=L\left(u\bigr|_{[0,1]}\right)=\int_0^1\Vert\dot u\Vert dx\ ,
\end{equation}
where $\Vert\cdot\Vert$ denotes the Weil-Petersson norm on tangent vectors to $\T$.  We use a similar convention for the energy $E(u)$:
\begin{equation} \label{E:energy}
E(u)=E\left(u\bigr|_{[0,1]}\right)=\int_0^1\Vert\dot u\Vert^2 dx\ .
\end{equation}
In particular, $L(u)\leq E^{1/2}(u)$.

 We define the  heat flow of $u$ as follows (cf.\ \cite{ES,Ham}). 
The flow is a time dependent family of
$\gamma$-equivariant paths $u(t,\cdot):\RBbb\to\T$, where for $t\geq 0$, $u(t,x)$ satisfies the heat equation:
\begin{equation}  \label{E:heatflow}
\dot u=\Delta u\ , \qquad u(0,x)=u(x)\ .
\end{equation}
Here, $\dot u$ denotes the $t$ derivative of $u(t,x)$, and $\Delta$ is the induced Laplacian from the Levi-Civit\`a connection of the
Weil-Petersson metric on
$\T$.  The following is standard:

\begin{Lem}  \label{L:heatflow}
For any initial condition, there is $T>0$ such that the solution $u(t,x)$ of (\ref{E:heatflow}) exists and is unique for all $0\leq t<T$.
\end{Lem}

To obtain existence for all time, we need to guarantee that the solution does not run off to the boundary $\partial \Tbar$ in finite time.  By
convexity, it suffices to require this at a single point, as the following lemma demonstrates:

\begin{Lem}  \label{L:longtime}
Let $u(t,x)$ be a $\gamma$-equivariant solution to the heat equation for $0\leq t<T$.
Suppose there is $x_0\in\RBbb$ and a compact set  $K\subset\T$ such that $u(t,x_0)\in K$  for all  $t< T$.  Then there is $\varepsilon>0$ such
that
$u(t,x)$ may be extended to a solution for $0\leq t<T+\varepsilon$.
\end{Lem}

\begin{proof}
Let $f$ be the convex exhaustion function on $\T$ from Proposition \ref{P:wolpert} (4).  Then by a straightforward calculation using the
negative curvature of $\T$,
$g(t,x)=f(u(t,x))$ is a subsolution to the heat equation.  By equivariance of the path we conclude that 
for each $N\in \ZBbb$
there is a compact set $K_N\subset\T$ such that $u(t,x_0+N)$ and $u(t,x_0-N)$ are contained in $K_N$ for all $t< T$.  In particular, 
 $g(t,x_0+N)$ and $g(t,x_0-N)$ are uniformly bounded in $t$.  On the other hand, since $g(t,x)$ is a
subsolution, an interior maximum on $[x_0-N,x_0+N]$ decreases in $t$.  Hence, for each fixed $N$, $g(t,x)$ is uniformly bounded for all
$0\leq t<T$ and 
 for all
$x\in [x_0-N,x_0+N]$.  So for $x$ in a compact set in $\RBbb$,
$u(t,x)$ lies in a compact set in $\T$.  It follows that as $t\to T$, $u(t,x)$ converges to a $\gamma$-equivariant path
$u_T:\RBbb\to\T$.  Using $u_T$ as initial conditions in Lemma \ref{L:heatflow}, we may find an extension to some
$T+\varepsilon$.
\end{proof}

We also have the following criterion for convergence of the flow at infinite time:

\begin{Lem}  \label{L:convergence}
Let $u(t,x)$ be a $\gamma$-equivariant solution to the heat equation.
Suppose there is $x_0\in\RBbb$, $T\geq 0$, and a compact set  $K\subset\T$ such that $u(t,x_0)\in K$  for all  $t\geq T$.  Then
$u(t,x)$ converges to a complete $\gamma$-equivariant geodesic in
$\T$.
\end{Lem}

\begin{proof}
The proof follows that of Lemma \ref{L:longtime} verbatim.  One obtains convergence of $u(t,x)$ along a subsequence $t_j\to \infty$  by standard
methods.  Convergence for all $t$ follows from Hartman's lemma \cite{Har}.
\end{proof}

Finally, given the existence  Theorem \ref{T:equivariantgeodesic}, we can state a general result for convergence of the heat flow:

\begin{Thm}  \label{T:convergence}
Let $\gamma\in\mcg$ be pseudo-Anosov.  There is a constant $c(\gamma) >0$ with the following significance:  If $u:\RBbb\to\T$ is a
$\gamma$-equivariant smooth map satisfying $E^{1/2}(u)\leq \twp(\gamma)+c(\gamma)$, and $u(t,x)$ denotes the heat flow (\ref{E:heatflow}) with
initial condition  $u(x)$, then $u(t,x)$ exists for all $t\geq 0$  and converges as $t\to \infty$ to the unique $\gamma$-equivariant complete
geodesic in
$\T$.
\end{Thm}

\begin{proof}   Let $A_\gamma:\RBbb\to\T$ denote the complete $\gamma$-equivariant geodesic from Theorem \ref{T:equivariantgeodesic}.  
By (\ref{E:distance}) it follows that if the length of a $\gamma$-equivariant map $u$ is near to $\twp(\gamma)$, then
$u$ must be uniformly close to $A_\gamma$.  We can now apply Hartman's distance decreasing property \cite{Har} to
conclude that along the flow, 
$$
\dwp(u(t,x), A_\gamma(x))\leq \sup_{x'}\,\dwp(u(x'), A_\gamma(x')) 
$$
for all $x$ and all $t$ in the domain of definition.  In particular, if the right hand side in the above expression is sufficiently
small compared to the distance of
$A_\gamma$ to
$\partial\Tbar$, then the inequality guarantees that $u(t,x)$ stays in a compact set in $\T$.  The result then follows from Lemmas
\ref{L:longtime} and \ref{L:convergence}.
\end{proof}

To conclude,
it is worth pointing out two additional consequences of the results in this paper.  First, we immediately obtain the following well-known fact concerning
mapping class groups (cf.\ \cite{Mc}):

\begin{Cor}  \label{C:centralizer}
Let $\gamma\in\mcg$ be pseudo-Anosov and let $\langle\gamma\rangle$ denote the cyclic group generated by $\gamma$.  Then the centralizer and
normalizer of $\langle\gamma\rangle$ in $\mcg$ are both virtually cyclic.
\end{Cor}

\begin{proof}
It clearly suffices to prove the statement for the normalizer, which we denote by $N\langle \gamma\rangle$. 
Recall from Section \ref{S:existence} that there is a finite index subgroup $\mcg^\prime\subset\mcg$ which contains no periodic mapping classes. 
Let $N'\langle\gamma\rangle=N\langle\gamma\rangle\cap\mcg^\prime$, and let $A_\gamma$ the axis of $\gamma$.  
Then for $g\in N^\prime\langle
\gamma\rangle$, we have $\gamma^k g A_\gamma\subset g A_\gamma$ for some $k$.  But $A_\gamma$ is also the unique $\gamma^k$-equivariant
geodesic, so $g A_\gamma=A_\gamma$.  In this way,  we have a homomorphism $N^\prime\langle \gamma\rangle\to \Iso(A_\gamma)\simeq \Iso(\RBbb)$,
which, by the assumption of no periodic mapping classes, must be injective.  The result follows, since by proper discontinuity of the action of
$\mcg$, the image of
$N'\langle\gamma\rangle$ in
$\Iso(A_\gamma)$ is also discrete.
\end{proof}

\noindent
Second, the proof of Theorem \ref{T:classification} shows that, just like for the Teichm\"uller metric, there is a
\emph{shortest} length for a closed Weil-Petersson geodesic in the moduli space of curves.  

\begin{Thm}  \label{T:shortest}
There exists a constant $\delta(g,n)>0$ depending only on $(g,n)$ such that if $\gamma\in\mcg$ satisfies $\twp(\gamma)\neq 0$, then
$\twp(\gamma)\geq \delta(g,n)$.
\end{Thm}

\begin{proof}
The proof proceeds by induction on $3g-3+n$.  Reducible elements fall into the inductive hypothesis.  Hence, we may assume that $\gamma$
is pseudo-Anosov.  
Let $A_\gamma$ be the axis of $\gamma$ as above. We claim that there is some $\delta_1>0$ depending only on $(g,n)$ such that either
$\twp(\gamma)\geq\delta_1$ or there is some point $\sigma\in A_\gamma$ such that $\distwp(\sigma,\partial\T)\geq \delta_1$.  This follows from Lemma
\ref{L:discontinuous} and the compactness of $\Mbar$.  Next, on the set of points a distance $\geq \delta_1/2$ from $\partial\T$, the Weil-Petersson and
Teichm\"uller metrics are uniformly quasi-isometric.  If $\twp(\gamma)\leq \delta_1/2$ is very small, then the translation length of $\gamma$ with
respect to the Teichm\"uller metric would also be very small.  On the other hand, the latter has a uniform lower bound depending only on $(g,n)$
(cf.\ \cite{Pe}).  This completes the proof.
\end{proof}


\section{Appendix}       \label{S:appendix}


The purpose of this appendix is to prove the asymptotic estimate for the Weil-Petersson metric asserted in (\ref{E:wp}).  The idea, due to Yamada
\cite{Y1}, is to combine the original construction of Masur \cite{M2} with the precise estimates for the hyperbolic metric due to Wolpert \cite{W5}.  This
requires some background notation which we shall sketch rather quickly.  A more detailed discussion may be found in the references cited.

Consider a point $\sigma\in \widehat{\bf T}={\bf T}_{g_1,n_1}\times\cdots\times{\bf T}_{g_N,n_N}$ obtained by collapsing curves $c_1,\ldots, c_k$. 
Then $\sigma$ corresponds to a nodal surface:  that is, a collection of Riemann surface $\Sigma_{g_1,n_1}, \ldots, \Sigma_{g_N,n_N}$ with
identifications that are holomorphically the coordinate lines in $\CBbb^2$ made in a neighborhood of certain of the marked points.  More precisely,
each curve $c_\ell$ gives rise to a pair of points $p_\ell, q_\ell$ that are identified with  marked points in the disjoint union of
$\Sigma_{g_j,n_j}$.  Let
$z_\ell, w_\ell$ indicate local conformal coordinates at $p_\ell, q_\ell$.  Then the (connected) nodal surface $\Sigma_0$ is formed via the
identification
$z_\ell w_\ell=0$.  We construct a degenerating family of Riemann surfaces $\Sigma_t$, $t=(t_1,\ldots, t_k)$, by replacing the nodal neighborhoods
with annuli $A_{\ell,t}$ given by the equations
$z_\ell w_\ell=t_\ell$.  We also may deform the conformal structure on the  initial nodal surface $\Sigma_0$ in regions compactly supported
away from the nodes.  In the notation of (\ref{E:wp}), this gives rise to deformations $\Sigma_{t,\tau}$ of the family $\Sigma_t$
corresponding to the coordinates
$\tau_{k+1},\ldots,
\tau_{3g-3+n}$.  Passing to the infinite abelian covering with coordinates
$(\theta_\ell,|t_\ell|)$, these degenerating families parametrize a neighborhood of $\sigma$ in the model neighborhood.  For brevity we
will omit the $\tau$ parameters from the notation.  Since the uniformizing coordinates may be chosen to vary continuously with respect to the $\tau$
parameters, this is no loss of generality.

We will denote by
$\Sigma^\ast$ the surface with boundary which is the complement of a neighborhood of the nodes of
$\Sigma_0$.  The gluing of $\Sigma^\ast$ with the annuli $A_{\ell,t}$ takes place on pairs of annuli $B_{\ell,t}$, also supported away from the nodes
(see \cite{W5}).
The nodal surface $\Sigma_0$ carries a complete hyperbolic metric $ds_0^2=\rho_0(z)|dz|^2$ which induces a hyperbolic metric on the surface
$\Sigma^\ast$ with boundary. We assume that the conformal coordinates used in the gluing have been chosen so that in a neighborhood of each node,
$\rho_0(z_\ell)=(|z_\ell|\log|z_\ell|)^{-2}$.
 Furthermore, each annulus $A_{\ell,t}$ carries a hyperbolic metric:
\begin{equation} \label{E:annulus}
ds^2_{\ell,|t|}= ds^2_{\ell,0} \frac{\Theta_{\ell,|t|}^2}{\left(\sin\Theta_{\ell,|t|}\right)^2}\ ,
\end{equation}
where
$$
ds_{\ell,0}^2=\frac{|dz_\ell|^2}{|z_\ell|^2(\log|z_\ell|)^2}\quad ,\quad \Theta_{\ell,|t|}=\pi\frac{\log|z_\ell|}{\log|t_\ell|}\ .
$$
Following Wolpert, we construct a ``grafted" metric $\graft$ on $\Sigma_t$ as follows:  $\graft$ agrees with the hyperbolic metrics $ds^2_{\ell,|t|}$
on each annulus $A_{\ell,t}$, $\graft$ is a smooth family of hyperbolic metrics with boundary on $\Sigma^\ast$ converging as $t\to 0$ to the one
induced by
$ds_0^2$, and in the regions
$B_{\ell,t}$,
$\graft$ interpolates between these hyperbolic metrics. The actual hyperbolic metric $ds^2_t$ on $\Sigma_t$ will be related to $\graft$ by a conformal factor:
$$
\rho_t(z)|dz|^2=ds_t^2=\frac{\rho_t(z)}{\rho_{t, \text{graft}}(z)}\graft\ .
$$
Then we have the following
uniform estimate on $\Sigma_t$ (cf.\ \cite[Expansion 4.2]{W5}):
\begin{equation} \label{E:w1}
\left| \frac{\rho_t(z)}{\rho_{t, \text{graft}}(z)}-1\right|=\sum_{\ell=1}^k O\left( (\log|t_\ell|)^{-2}\right)\ .
\end{equation}
We also note that in each $A_{\ell,t}$ the interpolation satisfies:
\begin{equation} \label{E:w2}
\left| \frac{\rho_0(z_\ell)}{\rho_{t, \text{graft}}(z_\ell)}-1\right|=O\left(\Theta_{\ell,|t|}^2 \right)
\end{equation}
(cf.\ \cite[\S 3.4.MG]{W5}).

Next, we review Masur's estimate of the Weil-Petersson metric.    According to \cite{M2}, we can find a family of quadratic differentials
$\varphi_i(z,t)$, meromorphic on $\Sigma_t$ with at most simple poles at the $n$ marked points, parametrizing the cotangent bundle of $\T$ at the
point represented by
$\Sigma_t$.   Then for $i\geq k+1$, the $\varphi_i$ converge as $t\to 0$ to quadratic differentials with at most simple
poles at the marked points (which now include the nodes).  For $i\leq k$, the $\varphi_i$ represent normal directions to the boundary stratum.

In terms of the local coordinate $z_\ell$   for the annular region $A_{\ell,t}$
 we have the following uniform bounds (cf.\ \cite{M2}, 5.4, 5.5, and note the remark on p.\ 633 before ``Proof of Theorem 1").  Here, $C$
denotes a constant that is independent of $t$.
\begin{align}
& \text{For}\ i,j\geq k+1\ ,\ |\varphi_i(z_\ell,t)\bar\varphi_j(z_\ell,t)|\leq C/|z_\ell|^2\label{E:m1} \\
&\text{For}\ i,j\leq k\ ,\
 |\varphi_i(z_\ell,t)\bar\varphi_j(z_\ell,t)|\leq
\begin{cases}
C_1|t_\ell|^2/|z_\ell|^4 + C|t_\ell|^2/|z_\ell|^3 &  i=j=\ell\ ,\ C_1\neq 0   \\
C|t_i||t_j|/|z_\ell|^3 &  i=\ell, j\neq \ell \ \text{or \emph{vice versa}} \\
C|t_i||t_j|/|z_\ell|^2 &  i\neq \ell, j\neq \ell
\end{cases}\label{E:m2}
\\
&\text{For}\ i\leq k\ , \ j\geq k+1\ ,\
 |\varphi_i(z_\ell,t)\bar\varphi_j(z_\ell,t)|\leq
\begin{cases}
C|t_i|/|z_\ell|^3 &  i=\ell \\
C|t_i|/|z_\ell|^2  &  i\neq \ell
\end{cases}\label{E:m3}
\\
&\text{For}\ i,j\geq k+1\ ,\
\left|\varphi_i(z_\ell,t)\bar\varphi_j(z_\ell,t)-\varphi_i(z_\ell,0)\bar\varphi_j(z_\ell,0)\right|\leq C|t_\ell|/|z_\ell|^4\label{E:m4}
\end{align}

Let $G_{i\bar j}(t)$ denote the components of the Weil-Petersson metric tensor with respect to this choice of (dual) basis.  
With this understood, we now state
\begin{Lem}  \label{L:metricestimate}
\begin{enumerate}
\item[(i)]  There is a nonzero constant $C$ independent of $t$ such that for $i\leq k$,
$$
\left| G_{i\bar i}(t)\right|= C|t_i|^{-2}(-\log|t_i|)^{-3}\left(1+\sum_{\ell=1}^k O\left((-\log|t_\ell|)^{-2}\right)\right)\ .
$$
\item[(ii)]  For $i,j\leq k$, $i\neq j$,
$$
\left| G_{i\bar j}(t)\right|= O\left(|t_i|^{-1}|t_j|^{-1}(-\log|t_i|)^{-3}(-\log|t_j|)^{-3}\right)\ .
$$
\item[(iii)]  For $i\leq k$, $j\geq k+1$,
$$
\left| G_{i\bar j}(t)\right|= O\left(|t_i|^{-1}(-\log|t_i|)^{-3}\right)\ .
$$
\item[(iv)]  For $i,j\geq k+1$,
$$
\left|G_{i\bar j}(t)-G_{i\bar j}(0)\right|=\sum_{\ell=1}^k O\left((-\log|t_\ell|)^{-2}\right)\ .
$$
\end{enumerate}
\end{Lem}

\noindent  Let us first point out that the estimate (\ref{E:wp}) follows from Lemma \ref{L:metricestimate} by making the substitution:
$(\theta_i,|t_i|)\mapsto (\theta_i,\xi_i)$, where $\xi_i=(-\log|t_i|)^{-1/2}$.  The rest of this section is therefore devoted to the 

\begin{proof}[Proof of Lemma \ref{L:metricestimate}]
As in \cite{M2}, the result will follow by first estimating the cometric.  Let $\langle\varphi_i(t),\varphi_j(t)\rangle$ denote the Weil-Petersson
pairing.  We will show that:
\begin{enumerate}
\item[(i')]   There is a nonzero constant $C$ such that for $i\leq k$,
$$
\left| \langle\varphi_i(t),\varphi_i(t)\rangle\right|= C|t_i|^{2}(-\log|t_i|)^{3}\left(1+\sum_{\ell=1}^k O\left((-\log|t_\ell|)^{-2}\right)\right)\ .
$$
\item[(ii')]   For $i,j\leq k$, $i\neq j$,
$
\left| \langle\varphi_i(t),\varphi_j(t)\rangle\right|=O\left(|t_i||t_j|\right)
$.
\item[(iii')]  For $i\leq k$, $j\geq k+1$,
$
\left| \langle\varphi_i(t),\varphi_j(t)\rangle\right|=O\left(|t_i|\right)
$.
\item[(iv')]  For $i,j\geq k+1$,
$$
\left|\langle\varphi_i(t),\varphi_j(t)\rangle-\langle\varphi_i(0),\varphi_j(0)\rangle\right|=\sum_{\ell=1}^k O\left((-\log|t_\ell|)^{-2}\right)\ .
$$
\end{enumerate}
Given these estimates, Lemma \ref{L:metricestimate} then follows via the cofactor expansion for the inverse  as in \cite{M2}.

To begin, notice that (ii') and (iii') are precisely the statements (iv) and (vi), respectively, of \cite[p.\ 634]{M2}.  To prove
(i'), we note that for $i\leq k$, $|\varphi_i(t)|^2$ is $O(|t_i|^2)$ uniformly on compact sets away from the pinching region.  We therefore need
to estimate:
\begin{equation}  \label{E:integral1}
\int_{A_{\ell,t}}\frac{|\varphi_i(z_\ell,t)|^2}{\rho_t(z_\ell)} |dz_\ell|^2\ .
\end{equation}
Using (\ref{E:w1})  we observe that (\ref{E:integral1}) is bounded by 
\begin{align}
\int_{A_{\ell,t}}\frac{|\varphi_i(z_\ell,t)|^2}{\rho_0(z_\ell)} |dz_\ell|^2
&+\int_{A_{\ell,t}}\frac{|\varphi_i(z_\ell,t)|^2}{\rho_0(z_\ell)}\left(\frac{\rho_0(z_\ell)}{\rho_t(z_\ell,t)}-1\right)|dz_\ell|^2
\notag \\
& \leq\
\int_{A_{\ell,t}}|\varphi_i(z_\ell,t)|^2|z_\ell|^2(\log|z_\ell|)^2|dz_\ell|^2\notag \\
&\qquad +C\int_{A_{\ell,t}}|\varphi_i(z_\ell,t)|^2|z_\ell|^2(\log|z_\ell|)^2\left(\frac{\rho_{t,
\text{graft}}(z_\ell)}{\rho_t(z_\ell,t)}-1\right)|dz_\ell|^2\notag \\
&\qquad +\int_{A_{\ell,t}}|\varphi_i(z_\ell,t)|^2|z_\ell|^2(\log|z_\ell|)^2\left(\frac{\rho_0(z_\ell)}{\rho_{t,
\text{graft}}(z_\ell,t)}-1\right)|dz_\ell|^2 \notag \\&
\leq\
\int_{A_{\ell,t}}|\varphi_i(z_\ell,t)|^2|z_\ell|^2(\log|z_\ell|)^2|dz_\ell|^2
 \label{E:integral3}
\\
&\qquad + \sum_{j=1}^k O\left( (\log|t_j|)^{-2}\right)
\int_{A_{l,t}}|\varphi_i(z_\ell,t)|^2|z_\ell|^2(\log|z_\ell|)^2|dz_\ell|^2 \notag\\
&+\int_{A_{\ell,t}}|\varphi_i(z_\ell,t)|^2|z_\ell|^2(\log|z_\ell|)^2\left(\frac{\rho_0(z_\ell)}{\rho_{t,
\text{graft}}(z_\ell,t)}-1\right)|dz_\ell|^2 \notag
\end{align}
Using (\ref{E:m2}) we see that the right hand side of (\ref{E:integral3}) is $O(|t_i|^2)$ if $i\neq l$. If
$i=l$, then again by (\ref{E:m2}) the first term on the right hand side of (\ref{E:integral3}) 
 is $O(|t_i|^2(-\log|t_i|)^3)$, and  this  contribution is in fact nonzero. The second term is similarly bounded.
Finally, the third term is bounded in the same way, after applying (\ref{E:w2}).
This proves (i').

To prove (iv'), first note that by (\ref{E:w1}) we have an estimate
$$
\left|
\int_K\frac{\varphi_i(z,t)\bar\varphi_j(z,t)}{\rho_t(z)} |dz|^2-\int_K\frac{\varphi_i(z,0)\bar\varphi_j(z,0)}{\rho_0(z)} |dz|^2
\right| = \sum_{\ell=1}^k O\left((\log|t_\ell|)^{-2}\right)
$$
where $K$ is a compact set supported away from the pinching region.  Note also that since $\varphi_i(z,0)$ have at most simple poles,
$$
\left|
\int_{|z_\ell| \leq |t_\ell|}\frac{\varphi_i(z_\ell,0)\bar\varphi_j(z_\ell,0)}{\rho_0(z_\ell)} |dz_\ell|^2
\right|= O(|t_\ell|)\ .
$$
 Hence, we need only estimate
\begin{align}
\biggl|
\int_{A_{\ell,t}}\frac{\varphi_i(z_\ell,t)\bar\varphi_j(z_\ell,t)}{\rho_t(z_\ell)} |dz_\ell|^2
&-\int_{A_{\ell,t}}\frac{\varphi_i(z_\ell,0)\bar\varphi_j(z_\ell,0)}{\rho_0(z_\ell)} |dz_\ell|^2
\biggr|
\notag \\
& \leq\
\int_{A_{\ell,t}}\frac{\left|\varphi_i(z_\ell,t)\bar\varphi_j(z_\ell,t)\right|}{\rho_0(z_\ell)}
\left|\frac{\rho_0(z_\ell)}{\rho_t(z_\ell,t)}-1\right||dz_\ell|^2
 \notag
\\
&\qquad +
\int_{A_{\ell,t}}\left|\varphi_i(z_\ell,t)\bar\varphi_j(z_\ell,t)-\varphi_i(z_\ell,0)\bar\varphi_j(z_\ell,0)\right|\frac{|dz_\ell|^2}
{\rho_0(z_\ell)} \notag \\
& \leq\
C\int_{A_{\ell,t}}\frac{\left|\varphi_i(z_\ell,t)\bar\varphi_j(z_\ell,t)\right|}{\rho_0(z_\ell)}
\left|\frac{\rho_{t,\text{graft}}(z_\ell)}{\rho_t(z_\ell,t)}-1\right||dz_\ell|^2
 \label{E:integral2}
\\
&\qquad +
\int_{A_{\ell,t}}\frac{\left|\varphi_i(z_\ell,t)\bar\varphi_j(z_\ell,t)\right|}{\rho_0(z_\ell)}
\left|\frac{\rho_0(z_\ell)}{\rho_{t,\text{graft}}(z_\ell,t)}-1\right||dz_\ell|^2
 \notag\\
&\qquad +
\int_{A_{\ell,t}}\left|\varphi_i(z_\ell,t)\bar\varphi_j(z_\ell,t)-\varphi_i(z_\ell,0)\bar\varphi_j(z_\ell,0)\right|\frac{|dz_\ell|^2}{\rho_0(z_\ell)}
\notag
\end{align}
By (\ref{E:w1}) and (\ref{E:m1}), the first term on the right hand side of (\ref{E:integral2}) is bounded by 
$$
\sum_{p=1}^k O\left((\log|t_p|)^{-2}\right)\int_{A_{\ell,t}}(\log|z_\ell|)^2
|dz_\ell|^2 = \sum_{p=1}^k O\left((\log|t_p|)^{-2}\right)\ .
$$ 
By (\ref{E:w2}), the second term is $O\left((\log|t_\ell|)^{-2}\right)$.
  Using (\ref{E:m4}), the third term on the right hand side of
(\ref{E:integral2}) is bounded by a constant times
$$
\int_{A_{\ell,t}}\frac{|t_\ell|}{|z_\ell|^4}|z_\ell|^2(\log|z_\ell|)^2|dz_\ell|^2 = O\left(|t_\ell|(-\log|t_\ell|)^{3}\right)\ .
$$
Hence, the desired estimate holds on each annular region  $A_{l,t}$ as well.  This proves (iv') and completes the proof of the lemma.
\end{proof}

\end{document}